\documentclass[11pt]{article}
\usepackage{mathrsfs}
\usepackage{amsfonts}
\usepackage{amssymb}
\usepackage{amsmath}
\usepackage[all]{xy}

\def\Ker{\mathop{\rm Ker}\nolimits}
\def\Hom{\mathop{\rm Hom}\nolimits}
\def\Coker{\mathop{\rm Coker}\nolimits}
\def\Mod{\mathop{\rm Mod}\nolimits}
\def\mod{\mathop{\rm mod}\nolimits}
\def\pd{\mathop{\rm pd}\nolimits}
\def\id{\mathop{\rm id}\nolimits}
\def\fd{\mathop{\rm fd}\nolimits}
\def\Gpd{\mathop{\rm Gpd}\nolimits}
\def\Gid{\mathop{\rm Gid}\nolimits}
\def\Gfd{\mathop{\rm Gfd}\nolimits}
\def\GFid{\mathop{\rm GFid}\nolimits}
\def\SGfd{\mathop{\rm SGfd}\nolimits}
\def\Im{\mathop{\rm Im}\nolimits}
\def\Tr{\mathop{\rm Tr}\nolimits}
\def\Ext{\mathop{\rm Ext}\nolimits}

\def\Mod{\mathop{\rm Mod}\nolimits}
\def\mod{\mathop{\rm mod}\nolimits}
\def\res{\mathop{\rm res}\nolimits}
\def\cores{\mathop{\rm cores}\nolimits}
\def\dim{\mathop{\rm dim}\nolimits}
\def\codim{\mathop{\rm codim}\nolimits}

\def\inf{\mathop{\rm inf}\nolimits}

\def\FP{\mathop{\rm FP}\nolimits}

\title{\Large \bf Homological Dimensions Relative to Preresolving Subcategories
}
\author{Zhaoyong Huang\thanks{{\it E-mail address}: huangzy@nju.edu.cn}
\\ {\footnotesize {\it Department of Mathematics, Nanjing University,
Nanjing 210093, Jiangsu Province, P.R. China}}}
\date{}
\begin{document}
\baselineskip=18pt \maketitle

\begin{abstract}
We introduce relative preresolving subcategories and precoresolving
subcategories of an abelian category and define homological
dimensions and codimensions relative to these subcategories
respectively. We study the properties of these homological
dimensions and codimensions and unify some important properties
possessed by some known homological dimensions. Then we apply the
obtained properties to special subcategories and in particular to
module categories. Finally we propose some open questions and
conjectures, which are closely related to the generalized Nakayama
conjecture and the strong Nakayama conjecture.
\end{abstract}

\vspace{0.5cm}

\centerline{\bf 1. Introduction}

\vspace{0.2cm}

In classical homological theory, homological dimensions are
important and fundamental invariants and every homological dimension
of modules is defined relative to some certain subcategory of
modules. For example, projective, flat and injective dimensions of
modules are defined relative to the categories of projective, flat
and injective modules respectively. When projective, flat and
injective modules are generalized to Gorenstein projective,
Gorenstein flat and Gorenstein injective modules respectively in
relative homological theory, Gorenstein projective, Gorenstein flat
and Gorenstein injective dimensions emerge; and in particular, they
share many nice properties of projective, flat and injective
dimensions respectively (e.g. [AB, C, CFH, CI, DLM, EJ1, EJ2, EJL,
GD, GT, HI, H2, HuH, LHX, MD, SSW, Z]). Then a natural question is:
if two homological (co)dimensions relative to a category and its
subcategory are defined, what is the relation between these two
homological (co)dimensions?  The purpose of this paper is to study
this question. We introduce relative preresolving subcategories and
precoresolving subcategories and define homological dimensions and
codimensions relative to these subcategories respectively. Then we
study their properties and unify some important properties possessed
by some known homological dimensions.

This paper is organized as follows.

In Section 2, we give some terminology and some preliminary results;
in particular, we give the definition of homological (co)dimension
relative to a certain full and additive subcategory of an abelian
category.

In Section 3, we first give the definition of (pre)resolving
subcategories of an abelian category. Then we give some criteria for
computing and comparing homological dimensions relative to different
preresolving subcategories. Let $\mathscr{E}$ and $\mathscr{T}$ be
additive and full subcategories of an abelian category $\mathscr{A}$
such that $\mathscr{T}$ is $\mathscr{E}$-preresolving with an
$\mathscr{E}$-proper generator $\mathscr{C}$. Assume that $0\to M
\to T_1\to T_0 \to A \to 0$ is an exact sequence in $\mathscr{A}$
with both $T_0$ and $T_1$ objects in $\mathscr{T}$. Then there
exists an exact sequence $0\to M \to T \to C \to A \to 0$ in
$\mathscr{A}$ with $T$ an object in $\mathscr{T}$ and $C$ an object
in $\mathscr{C}$; and furthermore, if the former exact sequence is
$\Hom_{\mathscr{A}}(X,-)$-exact for some object $X$ in
$\mathscr{A}$, then so is the latter one. As applications of this
result, we get that an object in $\mathscr{A}$ is an
$n{\text-}\mathscr{C}$-cosyzygy if and only if it is an
$n{\text-}\mathscr{T}$-cosyzygy; and also get that the
$\mathscr{T}$-dimension of an object $A$ in $\mathscr{A}$ is at most
$n$ if and only if there exists an exact sequence $0\to K_n \to
C_{n-1}\to C_{n-2}\to \cdots \to C_0 \to A \to 0$ in $\mathscr{A}$
with all $C_i$ objects in $\mathscr{C}$ and $K_n$ an object in
$\mathscr{T}$. In addition, we give some sufficient conditions under
which the $\mathscr{T}$-dimension and the $\mathscr{C}$-dimension of
an object in $\mathscr{A}$ are identical.

Section 4 is completely dual to Section 3.

In Section 5, we apply the results in Sections 3 and 4 to special
subcategories and in particular to module categories. Some known
results are generalized. Finally we propose some questions and
conjectures concerning the obtained results, which are closely
related to the generalized Nakayama conjecture and the strong
Nakayama conjecture.

Throughout this paper, $\mathscr{A}$ is an abelian category and all
subcategories of $\mathscr{A}$ are full and additive.

\vspace{0.5cm}

\centerline{\bf 2. Preliminaries}

\vspace{0.2cm}

In this section, we give some terminology and some preliminary
results.

\vspace{0.2cm}

{\bf Definition 2.1.} ([Hu]) Let $\mathscr{C}$ be a subcategory of
$\mathscr{A}$ and $n\geq 0$.

(1) If there exists an exact sequence $0\to M \to C_{n-1} \to
C_{n-2} \to \cdots \to C_0 \to A \to 0$ in $\mathscr{A}$ with all
$C_i$ objects in $\mathscr{C}$, then $M$ is called an {\it
$n{\text-}\mathscr{C}$-syzygy object} (of $A$), and $A$ is called an
{\it $n{\text-}\mathscr{C}$-cosyzygy object} (of $M$); in this case,
we denote by $M=\Omega^n_{\mathscr{C}}(A)$ and
$A=\Omega^{-n}_{\mathscr{C}}(M)$. We denote by
$\Omega^n_{\mathscr{C}}(\mathscr{A})$ (resp.
$\Omega^{-n}_{\mathscr{C}}(\mathscr{A})$) the subcategory of
$\mathscr{A}$ consisting of $n{\text-}\mathscr{C}$-syzygy (resp.
$n{\text-}\mathscr{C}$-cosyzygy) objects.

(2) For an object $A$ in $\mathscr{A}$, the $\mathscr{C}$-dimension
(resp. $\mathscr{C}$-codimension), denoted by
$\mathscr{C}{\text-}\dim A$ (resp. $\mathscr{C}{\text-}\codim A$),
is defined as inf$\{n\geq 0\mid$ there exists an exact sequence $0
\to C_{n} \to \cdots \to C_{1} \to C_{0} \to A \to 0$ (resp. $0 \to
A \to C^{0} \to C^{1} \to\cdots \to C^{n} \to 0$) in $\mathscr{A}$
with all $C_i$ (resp. $C^i$) objects in $\mathscr{C}\}$. Set
$\mathscr{C}{\text-}\dim A$ (resp. $\mathscr{C}{\text-}\codim
A)=\infty$ if no such integer exists.

\vspace{0.2cm}

Let $\mathscr{C}$ be a subcategory of $\mathscr{A}$. We denote by
$\mathscr{C}^{\bot}=\{A$ is an object in $\mathscr{A}\mid
\Ext_{\mathscr{A}}^i(C,A)=0$ for any object $C$ in $\mathscr{C}$ and
$i\geq 1\}$ and $^{\bot}\mathscr{C}=\{A$ is an object in
$\mathscr{A}\mid \Ext_{\mathscr{A}}^i(A,C)=0$ for any object $C$ in
$\mathscr{C}$ and $i\geq 1\}$.

\vspace{0.2cm}

{\bf Lemma 2.2.} {\it Let $\mathscr{C}$ and $\mathscr{D}$ be
subcategories of $\mathscr{A}$, and let $M$ be an object in
$^{\bot}\mathscr{C}$ and $M^{'}$ an object in
$\Omega_{\mathscr{C}}^{-n}(\mathscr{A})$ such that some
$\Omega_{\mathscr{C}}^{n}(M^{'})$ is an object in
$\mathscr{D}^{\bot}$. If $\mathscr{D}{\text-}\dim M\leq n(<\infty)$,
then $\Ext_{\mathscr{A}}^{i}(M,M^{'})=0$ for any $i\geq 1$.}

\vspace{0.2cm}

{\it Proof.} By assumption, there exists an exact sequence:
$$0\to M^{''} \to C_{n-1} \to \cdots \to C_{1}\to C_{0} \to M^{'}\to
0$$ in $\mathscr{A}$ with all $C_i$ objects in $\mathscr{C}$ and
$M^{''}$ an object in $\mathscr{D}^{\bot}$. Let $M$ be an object in
$^{\bot}\mathscr{C}$. Then $\Ext_{\mathscr{A}}^i(M,M^{'})\cong
\Ext_{\mathscr{A}}^{n+i}(M,M^{''})$ for any $i\geq 1$. If
$\mathscr{D}{\text-}\dim M\leq n(<\infty)$, then there exists an
exact sequence:
$$0\to D_n \to \cdots \to D_{1} \to D_{0} \to M \to 0$$
in $\mathscr{A}$ with all $D_i$ objects in $\mathscr{D}$. So
$\Ext_{\mathscr{A}}^{n+i}(M,M^{''}) \cong
\Ext_{\mathscr{A}}^i(D_n,M^{''})=0$ for any $i\geq 1$ and hence
$\Ext_{\mathscr{A}}^{i}(M,M^{'})=0$ for any $i\geq 1$. \hfill
$\square$

\vspace{0.2cm}

Let $\mathscr{E}$ be a subcategory of $\mathscr{A}$. Recall from
[EJ2] that a sequence:
$$\mathbb{S}: \cdots \to S_1 \to S_2 \to S_3 \to \cdots$$
in $\mathscr{A}$ is called {\it
$\Hom_{\mathscr{A}}(\mathscr{E},-)$-exact} (resp. {\it
$\Hom_{\mathscr{A}}(-,\mathscr{E})$-exact}) if {\it
$\Hom_{\mathscr{A}}(E,\mathbb{S})$} (resp.
$\Hom_{\mathscr{A}}(\mathbb{S},E)$) is exact for any object $E$ in
$\mathscr{E}$. An epimorphism (resp. a monomorphism) $f$ in
$\mathscr{A}$ is called {\it $\mathscr{E}$-proper} (resp. {\it
$\mathscr{E}$-coproper}) if it is
$\Hom_{\mathscr{A}}(\mathscr{E},-)$-exact (resp.
$\Hom_{\mathscr{A}}(-,\mathscr{E})$-exact).

\vspace{0.2cm}

{\bf Proposition 2.3.} {\it Let $\mathscr{C}$ and $\mathscr{E}$ be
subcategories of $\mathscr{A}$ and let $\mathscr{C}$ be closed under
kernels of ($\mathscr{E}$-proper) epimorphisms. If
$$0\to A_1\to A_2\to A_3 \to 0 \eqno{(2.1)}$$ is a
($\Hom_{\mathscr{A}}(\mathscr{E},-)$-exact) exact sequence in
$\mathscr{A}$ with $A_3$ an object in $\mathscr{C}$, then
$\mathscr{C}{\text-}\dim A_1\leq \mathscr{C}{\text-}\dim A_2$.}

\vspace{0.2cm}

{\it Proof.} Let $\mathscr{C}{\text-}\dim A_2=n(<\infty)$ and
$$0\to C_n \to \cdots \to C_1 \to C_0 \to A_2 \to 0$$ be an exact
sequence in $\mathscr{A}$ with all $C_i$ objects in $\mathscr{C}$.
By [Hu, Theorem 3.2], there exist exact sequences:
$$0\to C_n \to \cdots \to C_1 \to C \to A_1 \to 0$$ and
$$0\to C \to C_0 \to A_3 \to 0 \eqno{(2.2)}$$ From the proof of [Hu, Theorem
3.2] we see that if (2.1) is
$\Hom_{\mathscr{A}}(\mathscr{E},-)$-exact, then so is (2.2). Because
$\mathscr{C}$ is closed under kernels of ($\mathscr{E}$-proper)
epimorphisms and $A_3$ is an object in $\mathscr{C}$ by assumption,
$C$ is an object in $\mathscr{C}$ and $\mathscr{C}{\text-}\dim
A_1\leq n$. \hfill $\square$

\vspace{0.2cm}

Let $\mathscr{C}$ be a subcategory of $\mathscr{A}$. We denote by
$\mathscr{C}\bot\mathscr{C}$ if $\Ext_{\mathscr{A}}^i(C_1,C_2)=0$
for any objects $C_1,C_2$ in $\mathscr{C}$ and $i\geq 1$, and denote
by $\mathscr{C}{\text-}\dim^{<\infty}$ (resp.
$\mathscr{C}{\text-}\codim^{<\infty}$) the subcategory of
$\mathscr{A}$ consisting of objects with finite
$\mathscr{C}$-dimension (resp. $\mathscr{C}$-codimension).

\vspace{0.2cm}

{\bf Lemma 2.4.} {\it Let $\mathscr{C}$ be a subcategory of
$\mathscr{A}$ such that $\mathscr{C}\bot\mathscr{C}$ and
$\mathscr{C}{\text-}\dim^{<\infty}$ is closed under direct summands,
and let $0\to K \to C \to A \to 0$ be an exact sequence in
$\mathscr{A}$ with $\mathscr{C}{\text-}\dim A<\infty$ and $C$ an
object in $\mathscr{C}$. If $K$ is an object in
$\mathscr{C}^{\bot}$, then $\mathscr{C}{\text-}\dim K<\infty$.}

\vspace{0.2cm}

{\it Proof.} Because $\mathscr{C}{\text-}\dim A<\infty$, there
exists an exact sequence:
$$0\to M \to C_0 \to A \to 0$$ in
$\mathscr{A}$ with $C_0$ an object in $\mathscr{C}$ and
$\mathscr{C}{\text-}\dim M<\infty$. Consider the following pull-back
diagram:
$$\xymatrix{
& & & 0 \ar[d] & 0 \ar[d] & \\
& & & M \ar[d] \ar@{=}[r] & M \ar[d] & \\
& 0 \ar[r] & K \ar@{=}[d] \ar[r] & N \ar[d] \ar[r] & C_0 \ar[d] \ar[r] & 0 \\
& 0 \ar[r] & K \ar[r] & C \ar[d] \ar[r] & A \ar[r]\ar[d] & 0 \\
& & & 0 & 0 & }$$ Because $\mathscr{C}\bot\mathscr{C}$ and
$\mathscr{C}{\text-}\dim M<\infty$, it is easy to get that
$M\in\mathscr{C}^{\bot}$ by dimension shifting. So the middle column
in the above diagram splits, and hence $\mathscr{C}{\text-}\dim
N\leq \mathscr{C}{\text-}\dim M<\infty$ by [Hu, Lemma 3.1]. Because
$K$ is an object in $\mathscr{C}^{\bot}$ by assumption, the middle
row in the above diagram also splits and $K$ is isomorphic to a
direct summand of $N$. Thus $\mathscr{C}{\text-}\dim K<\infty$.
\hfill$\square$

\vspace{0.2cm}

{\bf Definition 2.5.} Let $\mathscr{C}\subseteq\mathscr{T}$ be
subcategories of $\mathscr{A}$.

(1) (cf. [SSW]) $\mathscr{C}$ is called a {\it generator} (resp.
{\it cogenerator}) for $\mathscr{T}$ if for any object $T$ in
$\mathscr{T}$, there exists an exact sequence $0\to T^{'} \to C \to
T \to 0$ (resp. $0\to T \to C \to T^{'} \to 0$) in $\mathscr{T}$
with $C$ an object in $\mathscr{C}$.

(2) Let $\mathscr{E}$ be a subcategory of $\mathscr{A}$.
$\mathscr{C}$ is called an {\it $\mathscr{E}$-proper generator}
(resp. {\it $\mathscr{E}$-coproper cogenerator}) for $\mathscr{T}$
if for any object $T$ in $\mathscr{T}$, there exists a
$\Hom_{\mathscr{A}}(\mathscr{E},-)$ (resp.
$\Hom_{\mathscr{A}}(-,\mathscr{E})$)-exact exact sequence $0\to
T^{'}\to C \to T \to 0$ (resp. $0\to T\to C \to T^{'} \to 0$) in
$\mathscr{A}$ such that $C$ is an object in $\mathscr{C}$ and
$T^{'}$ is an object in $\mathscr{T}$.

\vspace{0.2cm}

{\bf Lemma 2.6.} {\it Let $\mathscr{C}\subseteq\mathscr{T}$ be
subcategories of $\mathscr{A}$ such that $\mathscr{C}$ is a
cogenerator for $\mathscr{T}$, and let $0\to A_1 \to A_2 \to A_3 \to
0$ be an exact sequence in $\mathscr{A}$ such that both $A_2$ and
$A_3$ are objects in $\mathscr{T}^{\bot}$. If $A_1$ is an object in
$\mathscr{C}^{\bot}$, then $A_1$ is an object in
$\mathscr{T}^{\bot}$.}

\vspace{0.2cm}

{\it Proof.} Let $0\to A_1 \to A_2 \to A_3 \to 0$ be an exact
sequence in $\mathscr{A}$ such that both $A_2$ and $A_3$ are objects
in $\mathscr{T}^{\bot}$. Then $\Ext_{\mathscr{A}}^i(T,A_1)=0$ for
any object $T$ in $\mathscr{T}$ and $i\geq 2$. Because $\mathscr{C}$
is a cogenerator for $\mathscr{T}$ by assumption, there exists an
exact sequence:
$$0\to T\to C \to T^{'} \to 0$$ in $\mathscr{A}$ with $C$ an object
in $\mathscr{C}$ and $T^{'}$ an object in $\mathscr{T}$, which
yields an exact sequence:
$$\Ext_{\mathscr{A}}^i(C,A_1)\to \Ext_{\mathscr{A}}^i(T,A_1)\to \Ext_{\mathscr{A}}^{i+1}(T^{'},A_1)$$
for any $i\geq 1$. Note that $\Ext_{\mathscr{A}}^{i+1}(T^{'},A_1)=0$
for any $i\geq 1$ by the above argument. So, if $A_1$ is an object
in $\mathscr{C}^{\bot}$, then $\Ext_{\mathscr{A}}^i(T,A_1)=0$ for
any $i\geq 1$ and $A_1$ is an object in $\mathscr{T}^{\bot}$.
\hfill$\square$

\vspace{0.2cm}

{\bf Lemma 2.7.} {\it Let $\mathscr{C}\subseteq\mathscr{T}$ be
subcategories of $\mathscr{A}$ such that $\mathscr{C}$ is a
cogenerator for $\mathscr{T}$ and $\mathscr{C}$ is closed under
direct summands. Then $\mathscr{T}\bigcap
\mathscr{T}^{\bot}\subseteq\mathscr{C}$.}

\vspace{0.2cm}

{\it Proof.} Let $T$ be an object in $\mathscr{T}\bigcap
\mathscr{T}^{\bot}$. Then there exists a split exact sequence:
$$0\to T\to C \to T^{'} \to 0$$ in $\mathscr{A}$ with $C$ an object
in $\mathscr{C}$ and $T^{'}$ an object in $\mathscr{T}$. So $T$ is
isomorphic to a direct summand of $C$. Because $\mathscr{C}$ is
closed under direct summands by assumption, $T$ is an object in
$\mathscr{C}$. \hfill$\square$

\vspace{0.2cm}

Sather-Wagstaff, Sharif and White introduced the Gorenstein category
$\mathcal{G}(\mathscr{C})$ as follows.

\vspace{0.2cm}

{\bf Definition 2.8.} ([SSW]) Let $\mathscr{C}$ be a subcategory of
$\mathscr{A}$. The {\it Gorenstein subcategory}
$\mathcal{G}(\mathscr{C})$ of $\mathscr{A}$ is defined as
$\mathcal{G}(\mathscr{C})=\{A$ is an object in $\mathscr{A}\mid$
there exists an exact sequence $\cdots \to C_1 \to C_0 \to C^0 \to
C^1 \to \cdots$ in $\mathscr{A}$ with all terms objects in
$\mathscr{C}$, which is both
$\Hom_{\mathscr{A}}(\mathscr{C},-)$-exact and
$\Hom_{\mathscr{A}}(-,\mathscr{C})$-exact, such that $A\cong
\Im(C_0\to C^0)\}$.

\vspace{0.2cm}

The Gorenstein category unifies the following notions: modules of
Gorenstein dimension zero ([AB]), Gorenstein projective modules,
Gorenstein injective modules ([EJ1]), $V$-Gorenstein projective
modules, $V$-Gorenstein injective modules ([EJL]),
$\mathcal{W}$-Gorenstein modules ([GD]), and so on (see [Hu] for the
details).

\vspace{0.5cm}

\centerline{\bf 3. Computation and Comparison of Homological
Dimensions}

\vspace{0.2cm}

In this section, we introduce the notion of (pre)resolving
subcategories of $\mathscr{A}$. Then we give some criteria for
computing and comparing homological dimensions relative to different
preresolving subcategories.

\vspace{0.2cm}

{\bf Definition 3.1.} Let $\mathscr{E}$ and $\mathscr{T}$ be
subcategories of $\mathscr{A}$. Then $\mathscr{T}$ is called {\it
$\mathscr{E}$-preresolving} in $\mathscr{A}$ if the following
conditions are satisfied.

(1) $\mathscr{T}$ admits an $\mathscr{E}$-proper generator.

(2) $\mathscr{T}$ is {\it closed under $\mathscr{E}$-proper
extensions}, that is, for any
$\Hom_{\mathscr{A}}(\mathscr{E},-)$-exact exact sequence $0\to
A_1\to A_2 \to A_3 \to 0$ in $\mathscr{A}$, if both $A_1$ and $A_3$
are objects in $\mathscr{T}$, then $A_2$ is also an object in
$\mathscr{T}$.

An $\mathscr{E}$-preresolving subcategory $\mathscr{T}$ of
$\mathscr{A}$ is called {\it $\mathscr{E}$-resolving} if the
following condition is satisfied.

(3) $\mathscr{T}$ is {\it closed under kernels of
$\mathscr{E}$-proper epimorphisms}, that is, for any
$\Hom_{\mathscr{A}}(\mathscr{E},-)$-exact exact sequence $0\to
A_1\to A_2 \to A_3 \to 0$ in $\mathscr{A}$, if both $A_2$ and $A_3$
are objects in $\mathscr{T}$, then $A_1$ is also an object in
$\mathscr{T}$.

\vspace{0.2cm}

The following list shows that the class of
$\mathscr{E}$-(pre)resolving subcategories is rather large.

\vspace{0.2cm}

{\bf Example 3.2.} (1) Let $\mathscr{A}$ admit enough projective
objects and $\mathscr{E}$ the subcategory of $\mathscr{A}$
consisting of projective objects. Then a subcategory of
$\mathscr{A}$ closed under $\mathscr{E}$-proper extensions is just a
subcategory of $\mathscr{A}$ closed under extensions. Furthermore,
if $\mathscr{C}=\mathscr{E}$ in the above definition, then an
$\mathscr{E}$-preresolving subcategory is just a subcategory which
contains all projective objects and is closed under extensions, and
an $\mathscr{E}$-resolving subcategory is just a projectively
resolving subcategory in the sense of [H2].

(2) Let $\mathscr{C}$ be a subcategory of $\mathscr{A}$ with
$\mathscr{C}\bot\mathscr{C}$. Then by [SSW, Corollary 4.5], the
Gorenstein subcategory $\mathcal{G}(\mathscr{C})$ of $\mathscr{A}$
is a $\mathscr{C}$-preresolving subcategory of $\mathscr{A}$ with a
$\mathscr{C}$-proper generator $\mathscr{C}$; furthermore, if
$\mathscr{C}$ is closed under kernels of epimorphisms, then
$\mathcal{G}(\mathscr{C})$ is a $\mathscr{C}$-resolving subcategory
of $\mathscr{A}$ by [SSW, Theorem 4.12(a)].

(3) Let $R$ be a ring, $\Mod R$ the category of left $R$-modules and
$\mathcal{P}(\Mod R)$ the subcategory of $\Mod R$ consisting of
projective modules. Recall from [EJ2] that a pair of subcategories
$(\mathscr{X},\mathscr{Y})$ of $\Mod R$ is called a {\it cotorsion
pair} if $\mathscr{X}=\{X \in \Mod R \mid \Ext_R^1(X,Y)=0$ for any
$Y \in \mathscr{Y}\}$ and $\mathscr{Y}=\{Y \in \Mod R \mid
\Ext_R^1(X,Y)=0$ for any $X \in \mathscr{X}\}$. If
$(\mathscr{X},\mathscr{Y})$ is a cotorsion pair in $\Mod R$, then
$\mathscr{X}$ is a $\mathcal{P}(\Mod R)$-preresolving subcategory of
$\Mod R$ with a $\mathcal{P}(\Mod R)$-proper generator
$\mathcal{P}(\Mod R)$ ([EJ2]).

(4) Let $R$ be a ring and $\mathcal{F}(\Mod R)$ the subcategory of
$\Mod R$ consisting of flat modules. Then by [Hu, Lemma 3.1 and
Theorem 3.2], it is not difficult to see that the subcategory of
$\Mod R$ consisting of strongly Gorenstein flat modules (see [DLM]
or Section 5 below for the definition) is an $\mathcal{F}(\Mod
R)$-resolving subcategory of $\Mod R$ with an $\mathcal{F}(\Mod
R)$-proper generator $\mathcal{P}(\Mod R)$.

(5) Let $R$ be a ring. Then, the subcategory of $\Mod R$ consisting
of the modules $A$ satisfying $\Ext_R^i(A,P)=0$ for any
$P\in\mathcal{P}(\Mod R)$ and $i\geq 1$, is a $\mathcal{P}(\Mod
R)$-resolving subcategory of $\Mod R$ with a $\mathcal{P}(\Mod
R)$-proper generator $\mathcal{P}(\Mod R)$. Let $R$ be a left
noetherian ring, $\mod R$ the category of finitely generated left
$R$-modules and $\mathcal{P}(\mod R)$ the subcategory of $\mod R$
consisting of projective modules. Then the subcategory of $\mod R$
consisting of the modules $A$ satisfying $\Ext_R^i(A,R)=0$ for any
$i\geq 1$ is a $\mathcal{P}(\mod R)$-resolving subcategory of $\mod
R$ with a $\mathcal{P}(\mod R)$-proper generator $\mathcal{P}(\mod
R)$.

\vspace{0.2cm}

Unless stated otherwise, in the rest of this section, we fix a
subcategory $\mathscr{E}$ of $\mathscr{A}$ and an
$\mathscr{E}$-preresolving subcategory $\mathscr{T}$ of
$\mathscr{A}$ admitting an $\mathscr{E}$-proper generator
$\mathscr{C}$. We will give some criteria for computing the
$\mathscr{T}$-dimension of a given object $A$ in $\mathscr{A}$, and
then compare it with the $\mathscr{C}$-dimension of $A$.

The following two propositions play a crucial role in this section.

\vspace{0.2cm}

{\bf Proposition 3.3.} {\it Let
$$0\to M \to T_1\buildrel {f} \over \longrightarrow T_0 \to A \to 0\eqno{(3.1)}$$
be an exact sequence in $\mathscr{A}$ with both $T_0$ and $T_1$
objects in $\mathscr{T}$. Then we have

(1) There exists an exact sequence:
$$0\to M\to T \to C\to A\to 0\eqno{(3.2)}$$ in
$\mathscr{A}$ with $T$ an object in $\mathscr{T}$ and $C$ an object
in $\mathscr{C}$.

(2) If (3.1) is $\Hom_{\mathscr{A}}(X,-)$-exact for some object $X$
in $\mathscr{A}$, then so is (3.2).}

\vspace{0.2cm}

{\it Proof}. (1) Let
$$0\to M \to T_1\buildrel {f} \over \longrightarrow T_0 \to A \to 0$$
be an exact sequence in $\mathscr{A}$ with both $T_0$ and $T_1$
objects in $\mathscr{T}$. Because there exists a
$\Hom_{\mathscr{A}}(\mathscr{E},-)$-exact exact sequence:
$$0\to T_0^{'}\to C \to T_0 \to 0$$
in $\mathscr{A}$ with $C$ an object in $\mathscr{C}$ and $T_0^{'}$
an object in $\mathscr{T}$, we have the following pull-back diagram:
$$\xymatrix{
& & 0 \ar[d] & 0 \ar[d] &  & \\
& & T_0^{'} \ar[d] \ar@{=}[r] & T_0^{'} \ar[d] & \\
& 0 \ar[r] & W\ar[d] \ar[r] & C \ar[d]  \ar[r] & A \ar@{=}[d] \ar[r] & 0 \\
& 0 \ar[r] & \Im f\ar[d]\ar[r] & T_0 \ar[d] \ar[r] & A \ar[r] & 0 \\
& & 0 & 0 & & }$$ Then consider the following pull-back diagram:
$$\xymatrix{
& & & 0 \ar[d] & 0 \ar[d] & \\
& & & T_0^{'} \ar[d] \ar@{=}[r] & T_0^{'} \ar[d] & \\
& 0 \ar[r] & M \ar@{=}[d] \ar[r] & T \ar[d] \ar[r] & W \ar[d] \ar[r] & 0 \\
& 0 \ar[r] & M \ar[r] & T_1 \ar[d] \ar[r] & \Im f \ar[r]\ar[d] & 0 \\
& & & 0 & 0 & }$$ Because the middle column in the first diagram is
$\Hom_{\mathscr{A}}(\mathscr{E},-)$-exact, the first column in the
first diagram (that is, the third column in the second diagram) and
the middle column in the second diagram are also
$\Hom_{\mathscr{A}}(\mathscr{E},-)$-exact by [Hu, Lemma 2.4(1)].
Because both $T_0^{'}$ and $T_1$ are objects in $\mathscr{T}$, $T$
is also an object in $\mathscr{T}$. Connecting the middle rows in
the above two diagrams we get the desired exact sequence.

(2) If (3.1) is $\Hom_{\mathscr{A}}(X,-)$-exact for some object $X$
in $\mathscr{A}$, then so are the third rows in the above two
diagrams. So the middle rows in the above two diagrams and (3.2) are
also $\Hom_{\mathscr{A}}(X,-)$-exact by [Hu, Lemma 2.4(1)]. \hfill
$\square$

\vspace{0.2cm}

As an application of Proposition 3.3, we get the following

\vspace{0.2cm}

{\bf Proposition 3.4.} {\it Let $n\geq 1$ and
$$0\to M\to T_{n-1}\to T_{n-2}\to \cdots \to T_0\to A\to 0$$
be an exact sequence in $\mathscr{A}$ with all $T_i$ objects in
$\mathscr{T}$. Then there exist an exact sequence:
$$0\to N\to C_{n-1}\to C_{n-2}\to\cdots\to C_0\to A\to 0$$
and a $\Hom_{\mathscr{A}}(\mathscr{E},-)$-exact exact sequence:
$$0\to T\to N\to M\to 0$$ in $\mathscr{A}$ with
all $C_i$ objects in $\mathscr{C}$ and $T$ an object in
$\mathscr{T}$. In particular, an object in $\mathscr{A}$ is an
$n{\text-}\mathscr{C}$-cosyzygy if and only if it is an
$n{\text-}\mathscr{T}$-cosyzygy.}

\vspace{0.2cm}

{\it Proof.} We proceed by induction on $n$. The case for $n=1$ has
been proved in the proof of Proposition 3.3. Now suppose that $n\geq
2$ and we have an exact sequence: $$0\to M\to T_{n-1}\to T_{n-2}\to
\cdots \to T_0\to A\to 0$$ in $\mathscr{A}$ with all $T_i$ objects
in $\mathscr{T}$. Put $K=\Ker(T_{1}\to T_{0})$. By Proposition 3.3,
we get an exact sequence:
$$0\to K\to T_{1}^{'}\to C_{0}\to A\to 0$$
in $\mathscr{A}$ with $T_{1}^{'}$ an object in $\mathscr{T}$ and
$C_{0}$ an object in $\mathscr{C}$. Put $A^{'}=\Im(T_{1}^{'}\to
C_{0})$. Then we get an exact sequence:
$$0\to M\to T_{n-1}\to T_{n-2}\to
\cdots \to T_{2}\to T_{1}^{'}\to A^{'}\to 0$$ in $\mathscr{A}$. Thus
the assertion follows from the induction hypothesis. \hfill
$\square$

\vspace{0.2cm}

The following corollary is an immediate consequence of Proposition
3.4.

\vspace{0.2cm}

{\bf Corollary 3.5.}  {\it Let $M$ be an object in $\mathscr{A}$ and
$n\geq 0$. If $\mathscr{T}{\text-}\codim M=n$. Then there exists a
$\Hom_{\mathscr{A}}(\mathscr{E},-)$-exact exact sequence $0\to T \to
N \to M \to 0$ in $\mathscr{A}$ with $\mathscr{C}{\text-}\codim
N\leq n$ and $T$ an object in $\mathscr{T}$.}

\vspace{0.2cm}

{\it Proof.} Let $M$ be an object in $\mathscr{A}$ with
$\mathscr{T}{\text-}\codim M=n$. Applying Proposition 3.4 with $A=0$
we get a $\Hom_{\mathscr{A}}(\mathscr{E},-)$-exact exact sequence
$0\to T \to N \to M \to 0$ in $\mathscr{A}$ with
$\mathscr{C}{\text-}\codim N\leq n$ and $T$ an object in
$\mathscr{T}$. \hfill$\square$

\vspace{0.2cm}

We give a criterion for computing the $\mathscr{T}$-dimension of an
object in $\mathscr{A}$ as follows.

\vspace{0.2cm}

{\bf Theorem 3.6.} {\it The following statements are equivalent for
any object $A$ in $\mathscr{A}$ and $n \geq 0$.

(1) $\mathscr{T}{\text-}\dim A\leq n$.

(2) There exists an exact sequence: $$0\to K_n\to C_{n-1}\to
C_{n-2}\to\cdots\to C_0\to A\to 0$$ in $\mathscr{A}$ with all
objects $C_i$ in $\mathscr{C}$ and $K_n$ an object in
$\mathscr{T}$.}

\vspace{0.2cm}

{\it Proof.} $(2)\Rightarrow (1)$ is trivial.

$(1)\Rightarrow (2)$ We proceed by induction on $n$. The case for
$n=0$ is trivial. If $n=1$, then there exists an exact sequence:
$$0\to T_1 \to T_0\to A \to 0$$
in $\mathscr{A}$ with both $T_0$ and $T_1$ objects in $\mathscr{T}$.
Applying Proposition 3.3 with $M=0$, we get an exact sequence:
$$0 \to K_1 \to C_0\to A\to 0$$
in $\mathscr{A}$ with $C_0$ an object in $\mathscr{C}$ and $K_1$ an
object in $\mathscr{T}$.

Now suppose $n\geq 2$. Then there exists an exact sequence: $$0\to
T_n\to T_{n-1}\to \cdots\to T_1\to T_0\to A\to 0$$ in $\mathscr{A}$
with all $T_i$ objects in $\mathscr{T}$. Put $M=\Ker(T_1\to T_0)$.
By Proposition 3.3, we get an exact sequence:
$$0\to M\to T_1^{'}\to C_0\to A\to 0$$ in $\mathscr{A}$ with
$C_0$ an object in $\mathscr{C}$ and $T_1^{'}$
an object in $\mathscr{T}$. Put $B=\Im(T_1^{'}\to C_0)$. Then we get
an exact sequence:
$$0\to T_n\to T_{n-1}\to \cdots\to T_1^{'}\to B\to 0.$$
By the induction hypothesis, we get the following exact sequence:
$$0\to K_n\to
C_{n-1}\to\cdots\to C_1 \to B \to 0$$ in $\mathscr{A}$ with all
$C_i$ objects in $\mathscr{C}$ and $K_n$ an object in $\mathscr{T}$.
Thus we get the desired exact sequence:
$$0\to K_n\to C_{n-1}\to C_{n-2}\to\cdots\to C_1\to C_0\to A\to 0.$$ \hfill$\square$

\vspace{0.2cm}

The following result gives a criterion for computing the
$\mathscr{T}$-codimension of an object in $\mathscr{A}$. To some
extent, the proof of this result is dual to that of Theorem 3.6, so
we omit it.

\vspace{0.2cm}

{\bf Theorem 3.7.} {\it The following statements are equivalent for
any object $M$ in $\mathscr{A}$ and $n \geq 0$.

(1) $\mathscr{T}{\text-}\codim M\leq n$.

(2) There exists an exact sequence: $$0\to M\to K^0\to
C^{1}\to\cdots\to C^{n-1} \to C^n \to 0$$ in $\mathscr{A}$ with
$K^0$ an object in $\mathscr{T}$ and all $C^i$ objects in
$\mathscr{C}$, that is, there exists an exact sequence: $$0\to M\to
T\to H \to 0$$ in $\mathscr{A}$ with $T$ an object in $\mathscr{T}$
and $\mathscr{C}{\text-}\codim H\leq n-1$.}

\vspace{0.2cm}

%{\it Proof.} $(2)\Rightarrow (1)$ is trivial.

%$(1)\Rightarrow (2)$ We proceed by induction on $n$. The case for
%$n=0$ is trivial. If $n=1$, then there exists an exact sequence:
%$$0\to M \to T^0 \to T^1\to 0$$
%in $\mathscr{A}$ with both $T^0$ and $T^1$ objects in $\mathscr{T}$.
%Applying Proposition 3.3 with $A=0$, we get an exact sequence:
%$$0\to M \to K^0 \to C^1\to 0$$
%in $\mathscr{A}$ with $K^0$ an object in $\mathscr{T}$ and $C^1$ an
%object in $\mathscr{C}$.

%Now suppose $n\geq 2$. Then there exists an exact sequence: $$0\to
%M\to T^0\to T^{1}\to\cdots\to T^{n-1} \to T^n \to 0$$ in
%$\mathscr{A}$ with all $T^i$ objects in $\mathscr{T}$. Put
%$N=\Im(T^0\to T^1)$. By the induction hypothesis, we get the
%following exact sequence:
%$$0\to N\to G^1\to
%C^{2}\to\cdots\to C^{n-1} \to C^n \to 0$$ in $\mathscr{A}$ with
%$G^1$ an object in $\mathscr{T}$ and all $C^i$ objects in
%$\mathscr{C}$. Put $A=\Im(G^1\to C^{2})$. By Proposition 3.3, we get
%an exact sequence:
%$$0\to M\to K^0\to C^{1}\to A\to 0$$ in $\mathscr{A}$ with $K^0$ an object
%in $\mathscr{T}$ and $C^1$ an object in $\mathscr{C}$. Thus we get
%the desired exact sequence:
%$$0\to M\to K^0\to C^{1}\to\cdots\to C^{n-1} \to C^n \to 0.$$ \hfill$\square$

%\vspace{0.2cm}

The following result gives a sufficient condition such that the
$n{\text-}\mathscr{C}$-syzygy of an object in $\mathscr{A}$ with
$\mathscr{T}$-dimension at most $n$ is in $\mathscr{T}$, in which
the first assertion generalizes [AB, Lemma 3.12].

\vspace{0.2cm}

{\bf Theorem 3.8.} {\it Let $\mathscr{T}$ be closed under kernels of
($\mathscr{E}$-proper) epimorphisms and
$\mathscr{T}\subseteq{\mathscr{C}^{\bot}}$, and let $n\geq 0$. Then
for any object $A$ in $\mathscr{A}$ with $\mathscr{T}{\text-}\dim
A\leq n$ we have

(1) For any ($\Hom_{\mathscr{A}}(\mathscr{E},-)$-exact) exact
sequence $0\to K_n \to C_{n-1}\to \cdots \to C_1\to C_0 \to A \to 0$
in $\mathscr{A}$ with all $C_i$ objects in $\mathscr{C}$, $K_n$ is
an object in $\mathscr{T}$.

(2) If $0\to K\to C \to A \to 0$ is a
($\Hom_{\mathscr{A}}(\mathscr{E},-)$-exact) exact sequence in
$\mathscr{A}$ with $C$ an object in $\mathscr{C}$, then
$\mathscr{T}{\text-}\dim K\leq n-1$.}

\vspace{0.2cm}

{\it Proof.} Let $\mathscr{T}{\text-}\dim A\leq n$ and
$\mathscr{T}\subseteq{\mathscr{C}^{\bot}}$. Then there exists a
$\Hom_{\mathscr{A}}(\mathscr{C},-)$-exact exact sequence:
$$0\to T_n \to C_{n-1}^{'}\to \cdots \to C_1^{'}\to C_0^{'} \to A \to 0$$
in $\mathscr{A}$ with all $C_i^{'}$ objects in $\mathscr{C}$ and
$T_n$ an object in $\mathscr{T}$ by Proposition 3.4.

(1) Let
$$0\to K_n \to C_{n-1}\to \cdots \to C_1\to C_0 \to A \to 0$$
be a ($\Hom_{\mathscr{A}}(\mathscr{E},-)$-exact) exact sequence in
$\mathscr{A}$ with all $C_i$ objects in $\mathscr{C}$. Then by [Hu,
Theorem 3.2] we get a ($\Hom_{\mathscr{A}}(\mathscr{E},-)$-exact)
exact sequence:
$$0\to K_n \to T_n\bigoplus C_{n-1}\to C_{n-1}^{'}\bigoplus C_{n-2}\cdots \to
C_{1}^{'}\bigoplus C_{0} \to C_{0}^{'} \to 0.$$ Because
$\mathscr{T}$ is closed under kernels of ($\mathscr{E}$-proper)
epimorphisms by assumption, $K_n$ is an object in $\mathscr{T}$.

(2) Put $T_1=\Im(C_1^{'}\to C_0^{'})$. Then we have a
$\Hom_{\mathscr{A}}(\mathscr{C},-)$-exact exact sequence:
$$0\to T_1\to C_0^{'} \to A \to 0$$ in $\mathscr{A}$ with $\mathscr{T}{\text-}\dim T_1\leq
n-1$. Let $0\to K\to C \to A \to 0$ be a
($\Hom_{\mathscr{A}}(\mathscr{E},-)$-exact) exact sequence in
$\mathscr{A}$ with $C$ an object in $\mathscr{C}$. Consider the
following pull-back diagram:
$$\xymatrix{
& & & 0 \ar[d] & 0 \ar[d] & \\
& & & T_1 \ar[d] \ar@{=}[r] & T_1 \ar[d] & \\
& 0 \ar[r] & K \ar@{=}[d] \ar[r] & X \ar[d] \ar[r] & C_0^{'} \ar[d] \ar[r] & 0 \\
& 0 \ar[r] & K \ar[r] & C \ar[d] \ar[r] & A \ar[r]\ar[d] & 0 \\
& & & 0 & 0 & }$$ Because the third column in this diagram is
$\Hom_{\mathscr{A}}(\mathscr{C},-)$-exact, the middle column is also
$\Hom_{\mathscr{A}}(\mathscr{C},-)$-exact by [Hu, Lemma 2.4(1)]. So
the middle column splits and $X\cong T_1 \bigoplus C$. Then the
middle row yields an exact sequence:
$$0\to K \to T_1\bigoplus C\to C_{0}^{'} \to 0.$$
By Proposition 2.3, $\mathscr{T}{\text-}\dim K\leq
\mathscr{T}{\text-}\dim T_1\bigoplus C\leq n-1$. \hfill$\square$

\vspace{0.2cm}

%The following corollary is an immediate consequence of Theorem 3.8,
%in which the first assertion in the following result generalizes
%[AB, Lemma 3.12].

%\vspace{0.2cm}

%{\bf Corollary 3.9.} {\it Let $\mathscr{T}$ admit a projective
%generator and be closed under extensions and kernels of
%epimorphisms. Then for any object $A$ in $\mathscr{A}$ and $n\geq 0$
%we have

%(1) $\mathscr{T}{\text-}\dim A\leq n$ if and only if for any exact
%sequence $0\to K_n \to C_{n-1}\to \cdots \to C_1\to C_0 \to A \to 0$
%in $\mathscr{A}$ with all $C_i$ projective objects in $\mathscr{T}$,
%$K_n$ is an object in $\mathscr{T}$.

%(2) If $0\to K\to P \to A \to 0$ is an exact sequence in
%$\mathscr{A}$ with $\mathscr{T}{\text-}\dim A\leq n$ and $P$ a
%projective object in $\mathscr{T}$, then $\mathscr{T}{\text-}\dim
%K\leq n-1$.}

%\vspace{0.2cm}

We use $\mathscr{T}{\text-}\dim^{\leq n}$ to denote the subcategory
of $\mathscr{A}$ consisting of objects with $\mathscr{T}$-dimension
at most $n$.

\vspace{0.2cm}

{\bf Corollary 3.9.} {\it Let $\mathscr{T}$ be a
$\mathscr{C}$-resolving subcategory of $\mathscr{A}$ with a
$\mathscr{C}$-proper generator $\mathscr{C}$ and
$\mathscr{T}\subseteq{\mathscr{C}^{\bot}}$. If $\mathscr{T}$ is
closed under direct summands, then so is
$\mathscr{T}{\text-}\dim^{\leq n}$ for any $n\geq 0$.}

\vspace{0.2cm}

{\it Proof.} The case for $n=0$ follows from the assumption. Now Let
$n\geq 1$ and let $A$ be an object in $\mathscr{A}$ with
$\mathscr{T}{\text-}\dim A\leq n$ and $A=A_1\bigoplus A_2$. Because
$\mathscr{T}\subseteq{\mathscr{C}^{\bot}}$ by assumption, there
exists a $\Hom_{\mathscr{A}}(\mathscr{C},-)$-exact exact sequence:
$$0\to K_n\to C_{n-1}\to
C_{n-2}\to\cdots\to C_0\buildrel {f_0}\over \to A\to 0$$ in
$\mathscr{A}$ with all $C_i$ objects in $\mathscr{C}$ and $K_n$ an
object in $\mathscr{T}$ by Theorem 3.6. Note that both $$0\to A_2
\buildrel {{\binom 0 {1_{A_2}}}} \over \longrightarrow A \buildrel
{(1_{A_1},0)} \over \longrightarrow A_1 \to 0$$ and $$0\to A_1
\buildrel {{\binom {1_{A_1}} 0}} \over \longrightarrow A \buildrel
{(0, 1_{A_2})} \over \longrightarrow A_2 \to 0$$ are exact and
split. So both
$$C_0 \buildrel {(1_{A_{1}},0) f_{0}} \over \longrightarrow A_1 \to 0$$
and
$$C_0 \buildrel {(0, 1_{A_{2}}) f_{0}} \over \longrightarrow A_2 \to 0$$ are
$\Hom_{\mathscr{A}}(\mathscr{C},-)$-exact exact sequences. By [Hu,
Theorem 3.6], we get the following
$\Hom_{\mathscr{A}}(\mathscr{C},-)$-exact exact sequences:
$$C_0\bigoplus C_1 \to C_0 \to A_1 \to 0$$
and
$$C_0\bigoplus C_1 \to C_0 \to A_2 \to 0.$$ Again by [Hu, Theorem
3.6], we get the following $\Hom_{\mathscr{A}}(\mathscr{C},-)$-exact
exact sequences:
$$C_0\bigoplus C_1\bigoplus C_2 \to C_0\bigoplus C_1 \to C_0 \to A_1 \to 0$$
and
$$C_0\bigoplus C_1\bigoplus C_2 \to C_0\bigoplus C_1 \to C_0 \to A_2 \to 0.$$
Continuing this procedure, we finally get the following
$\Hom_{\mathscr{A}}(\mathscr{C},-)$-exact exact sequences:
$$0\to X_n\to \bigoplus _{i=0}^{n-1}C_{i}\to
\bigoplus _{i=0}^{n-2}C_{i}\to\cdots\to C_1\bigoplus C_0 \to C_0\to
A_1\to 0$$ and
$$0\to Y_n\to \bigoplus _{i=0}^{n-1}C_{i}\to
\bigoplus _{i=0}^{n-2}C_{i}\to\cdots\to C_1\bigoplus C_0 \to C_0\to
A_2\to 0.$$ Put $U_j=\bigoplus _{i=0}^{j}C_{i}$ for any $0\leq j\leq
n-1$. Then we get a $\Hom_{\mathscr{A}}(\mathscr{C},-)$-exact exact
sequence:
$$0\to X_n\bigoplus Y_n\to U_{n-1}\bigoplus U_{n-1}\to
U_{n-2}\bigoplus U_{n-2}\to\cdots\to U_{1}\bigoplus U_{1} \to
U_{0}\bigoplus U_{0}\to A\to 0.$$ By Theorem 3.8, $X_n\bigoplus Y_n$
is an object in $\mathscr{T}$. So both $X_n$ and $Y_n$ are objects
in $\mathscr{T}$ and hence $\mathscr{T}{\text-}\dim A_1\leq n$ and
$\mathscr{T}{\text-}\dim A_2\leq n$. \hfill$\square$

\vspace{0.2cm}

The following result gives some sufficient conditions such that the
$\mathscr{T}$-dimension and the $\mathscr{C}$-dimension of an object
in $\mathscr{A}$ are identical.

\vspace{0.2cm}

{\bf Theorem 3.10.} {\it Let
$\mathscr{T}\subseteq{\mathscr{C}^{\bot}}\bigcap
{^{\bot}\mathscr{C}}$ and $\mathscr{C}$ be closed under direct
summands. Then for an object $A$ in $\mathscr{A}$,
$\mathscr{T}{\text-}\dim A=\mathscr{C}{\text-}\dim A$ if one of the
following conditions is satisfied.

(1) $\mathscr{C}{\text-}\dim A<\infty$, $\mathscr{E}=\mathscr{C}$
and $\mathscr{T}$ is closed under kernels of $\mathscr{C}$-proper
epimorphisms.

(2) $\mathscr{C}{\text-}\dim A<\infty$, $\mathscr{E}=\mathscr{C}$
and $\mathscr{C}{\text-}\dim^{<\infty}$ is closed under direct
summands.

(3) $A$ is an object in $\mathscr{T}^{\bot}$ and $\mathscr{C}$ is a
cogenerator for $\mathscr{T}$.}

\vspace{0.2cm}

{\it Proof.} It is trivial that $\mathscr{C}{\text-}\dim
A\geq\mathscr{T}{\text-}\dim A$. In the following we prove
$\mathscr{C}{\text-}\dim A\leq\mathscr{T}{\text-}\dim A$. Suppose
$\mathscr{T}{\text-}\dim A=n<\infty$.

(1) Let $\mathscr{C}{\text-}\dim A=t(<\infty)$. If $n<t$, then
consider the following $\Hom_{\mathscr{A}}(\mathscr{C},-)$-exact
exact sequence:
$$0\to C_t \to \cdots \to C_n\to C_{n-1}\to\cdots \to C_1 \to C_0 \to A \to 0$$
in $\mathscr{A}$ with all $C_i$ objects in $\mathscr{C}$. Put
$K_n=\Im(C_n\to C_{n-1})$. So $K_n$ is an object in $\mathscr{T}$ by
Theorem 3.8(1), and hence an object in ${^{\bot}\mathscr{C}}$ by
assumption. It yields that the exact sequence:
$$0\to C_t \to \cdots \to C_n\to K_n\to 0$$ splits and $K_n$ is
isomorphic to a direct summand of $C_n$. Because $\mathscr{C}$ is
closed under direct summands by assumption, $K_n$ is an object in
$\mathscr{C}$ and $\mathscr{C}{\text-}\dim A\leq n$, which is a
contradiction. So $n\geq t$.

In the following, we prove (2) and (3).

Let
$$0\to T_n \to T_{n-1}\to \cdots \to T_1 \to T_0 \to A \to 0$$
be an exact sequence in $\mathscr{A}$ with all $T_i$ objects in
$\mathscr{T}$. By Proposition 3.4, we get an exact sequence:
$$0\to K_n \to C_{n-1}\to \cdots \to C_1 \to C_0 \to A \to 0$$
and a $\Hom_{\mathscr{A}}(\mathscr{E},-)$-exact exact sequence:
$$0\to T \to K_n \to T_n \to 0$$
in $\mathscr{A}$ with all $C_i$ objects in $\mathscr{C}$ and $T$ an
object in $\mathscr{T}$. So $K_n$ is an object in
$\mathscr{T}(\subseteq{\mathscr{C}^{\bot}}\bigcap
{^{\bot}\mathscr{C}})$.

(2) Because $\mathscr{C}{\text-}\dim A<\infty$,
$\mathscr{C}{\text-}\dim K_n<\infty$ by Lemma 2.4. By assumption
$\mathscr{T}\subseteq{^{\bot}\mathscr{C}}$, it is easy to see that
$K_n$ is isomorphic to a direct summand of some object in
$\mathscr{C}$. Since $\mathscr{C}$ is closed under direct summands
by assumption, $K_n$ is an object in $\mathscr{C}$ and
$\mathscr{C}{\text-}\dim A\leq n$.

(3) Let $A$ be an object in $\mathscr{T}^{\bot}$ and $K_i=\Im(C_i
\to C_{i-1})$ for any $1\leq i \leq n-1$. Then all $K_i$ are objects
in $\mathscr{C}^{\bot}$. By Lemma 2.6, all $K_i$ are objects in
$\mathscr{T}^{\bot}$, and in particular $K_n$ is an object in
$\mathscr{T}^{\bot}$. So $K_n$ is an object in $\mathscr{C}$ by
Lemma 2.7, and hence $\mathscr{C}{\text-}\dim A\leq n$.
\hfill$\square$

\vspace{0.2cm}

The following result gives a sufficient condition such that the
$\mathscr{T}$-codimension and the $\mathscr{C}$-codimension of an
object in $\mathscr{A}$ are identical.

\vspace{0.2cm}

{\bf Theorem 3.11.} {\it Let $\mathscr{D}$ be a subcategory of
$\mathscr{A}$ such that
$\mathscr{T}\subseteq{^{\bot}\mathscr{C}}\bigcap{\mathscr{D}^{\bot}}$,
and let $\mathscr{C}{\text-}\codim^{\leq n}$ be closed under direct
summands for any $n\geq 0$. If $M$ is an object in $\mathscr{A}$
with $\mathscr{D}{\text-}\dim M<\infty$, then
$\mathscr{T}{\text-}\codim M=\mathscr{C}{\text-}\codim M$.}

\vspace{0.2cm}

{\it Proof.} It is clear that $\mathscr{C}{\text-}\codim M\geq
\mathscr{T}{\text-}\codim M$. In the following we prove
$\mathscr{T}{\text-}\codim M\geq\mathscr{C}{\text-}\codim M$.

Without loss of generality, assume $\mathscr{T}{\text-}\codim
M=n<\infty$. If $n=0$, then $M$ is an object in $\mathscr{T}$ and
there exists a $\Hom_{\mathscr{A}}(\mathscr{E},-)$-exact exact
sequence:
$$0\to M^{'}\to C \to M \to 0\eqno{(3.3)}$$ in $\mathscr{A}$ with $C$ an object
in $\mathscr{C}$ and $M^{'}$ an object in $\mathscr{T}$. Notice that
$\mathscr{T}\subseteq{^{\bot}\mathscr{C}\bigcap\mathscr{D}^{\bot}}$
by assumption, so $\Ext_{\mathscr{A}}^{i}(M,M^{'})=0$ for any $i\geq
1$ by Lemma 2.2. It follows that the exact sequence (3.3) splits,
which implies that $M$ is isomorphic to a direct summand of $C$.
Because $\mathscr{C}$ is closed under direct summands by assumption,
$M$ is an object in $\mathscr{C}$.

Now suppose $n\geq 1$. By Theorem 3.7, there exists an exact
sequence:
$$0\to M \to T \to H \to 0$$ in $\mathscr{A}$ with $T$ an object
in $\mathscr{T}$ and $\mathscr{C}{\text-}\codim H\leq n-1$. It is
easy to see that $M$ is an object in $^{\bot}\mathscr{C}$. Because
there exists a $\Hom_{\mathscr{A}}(\mathscr{E},-)$-exact exact
sequence:
$$0\to T^{'}\to C^{'} \to T \to 0$$ in $\mathscr{A}$ with $C^{'}$ an object
in $\mathscr{C}$ and $T^{'}$ an object in $\mathscr{T}$, we have the
following pull-back diagram:
$$\xymatrix{
& & 0 \ar[d] & 0 \ar[d] &  & \\
& & T^{'} \ar[d] \ar@{=}[r] & T^{'} \ar[d] & \\
& 0 \ar[r] & N\ar[d] \ar[r] & C^{'} \ar[d]  \ar[r] & H \ar@{=}[d] \ar[r] & 0 \\
& 0 \ar[r] & M\ar[d]\ar[r] & T \ar[d] \ar[r] & H \ar[r] & 0 \\
& & 0 & 0 & & }$$ By the exactness of the middle row in the above
diagram, $\mathscr{C}{\text-}\codim N\leq n$. Because
$\Ext_{\mathscr{A}}^{1}(M,T^{'})=0$ by Lemma 2.2, the first column
in the above diagram splits. So $M$ is isomorphic to a direct
summand of $N$. Because $\mathscr{C}{\text-}\codim^{\leq n}$ is
closed under direct summands by assumption,
$\mathscr{C}{\text-}\codim M\leq n$. \hfill$\square$

\vspace{0.2cm}

In the following, we fix a subcategory $\mathscr{C}$ of
$\mathscr{A}$.

The following two corollaries give some sufficient conditions such
that the $\mathcal{G}(\mathscr{C})$-dimension and the
$\mathscr{C}$-dimension of an object in $\mathscr{A}$ are identical.
The first one is a generalization of [Z, Theorem 2.3].

\vspace{0.2cm}

{\bf Corollary 3.12.} {\it Let $\mathscr{C}\bot\mathscr{C}$ and let
$\mathscr{C}$ be closed under direct summands. Then for any object
$A$ in ${\mathcal{G}(\mathscr{C})}^{\bot}$,
$\mathcal{G}(\mathscr{C}){\text-}\dim A=\mathscr{C}{\text-}\dim A$.}

\vspace{0.2cm}

{\it Proof.} Let $\mathscr{C}\bot\mathscr{C}$. It is clear that
$\mathscr{C}$ is a $\mathscr{C}$-proper generator and a
$\mathscr{C}$-coproper cogenerator for $\mathcal{G}(\mathscr{C})$.
By [SSW, Corollary 4.5], $\mathcal{G}(\mathscr{C})$ is closed under
extensions. By [Hu, Lemma 5.7], $\mathcal{G}(\mathscr{C})\subseteq
{\mathscr{C}^{\bot}}\bigcap {^{\bot}\mathscr{C}}$. Now the assertion
follows from Theorem 3.10(3). \hfill$\square$

\vspace{0.2cm}

The following is a generalization of [H1, Theorem 2.2] and [Z,
Corollary 2.5].

\vspace{0.2cm}

{\bf Corollary 3.13.} {\it Let $\mathscr{C}\bot\mathscr{C}$ and let
$\mathscr{C}$ be closed under direct summands, and let $\mathscr{D}$
be a subcategory of ${\mathcal{G}(\mathscr{C})}^{\bot}$. Then for
any object $A$ in $\mathscr{C}^{\bot}$ with
$\mathscr{D}{\text-}\codim A<\infty$,
$\mathcal{G}(\mathscr{C}){\text-}\dim A=\mathscr{C}{\text-}\dim A$.}

\vspace{0.2cm}

{\it Proof.} Let $A$ be an object in $\mathscr{C}^{\bot}$ with
$\mathscr{D}{\text-}\codim A<\infty$. Because $\mathscr{D}$ is a
subcategory of $\mathcal{G}(\mathscr{C})^{\bot}$, it is easy to see
that $A$ is an object in $\mathcal{G}(\mathscr{C})^{\bot}$ by Lemma
2.6. Then the assertion follows from Corollary 3.12. \hfill
$\square$

\vspace{0.2cm}

The following result gives a sufficient condition such that the
$\mathcal{G}(\mathscr{C})$-dimension and the
$^{\bot}\mathscr{C}$-dimension of an object in $\mathscr{A}$ are
identical.

\vspace{0.2cm}

{\bf Theorem 3.14.} {\it Let $\mathscr{C}\bot\mathscr{C}$ and let
$A$ an object in $\mathscr{A}$ with
$\mathcal{G}(\mathscr{C}){\text-}\dim A<\infty$. Then
$\mathcal{G}(\mathscr{C}){\text-}\dim
A={^{\bot}\mathscr{C}}{\text-}\dim A$.}

\vspace{0.2cm}

{\it Proof.} By [Hu, Lemma 5.7],
$\mathcal{G}(\mathscr{C}){\text-}\dim
A\geq{^{\bot}\mathscr{C}}{\text-}\dim A$. In the following we prove
$\mathcal{G}(\mathscr{C}){\text-}\dim A\leq
{^{\bot}\mathscr{C}}{\text-}\dim A$.

Suppose ${^{\bot}\mathscr{C}}{\text-}\dim A=n<\infty$ and
$\mathcal{G}(\mathscr{C}){\text-}\dim A=m<\infty$. If $n=0$, then
$A$ is an object in ${^{\bot}\mathscr{C}}$. So by [Hu, Theorem 5.8],
$A$ is an object in $\mathcal{G}(\mathscr{C})$ and $m=0$. Let $n\geq
1$. Then $\Ext_{\mathscr{A}}^{n+i}(A,C)=0$ for any object $C$ in
$\mathscr{C}$ and $i\geq 1$. If $m>n$, then consider the following
exact sequence:
$$0\to G_m \to \cdots \to G_n \to G_{n-1} \to G_1 \to G_0 \to A \to 0$$
in $\mathscr{A}$ with all $G_i$ objects in
$\mathcal{G}(\mathscr{C})$. Putting $K_n=\Im(G_n \to G_{n-1})$, then
$K_n$ is an object in ${^{\bot}\mathscr{C}}$ and
$\mathcal{G}(\mathscr{C}){\text-}\dim K_n\leq m-n<\infty$. By the
above argument, $K_n$ is an object in $\mathcal{G}(\mathscr{C})$. So
$\mathcal{G}(\mathscr{C}){\text-}\dim A\leq n$, which is a
contradiction. Thus $m\leq n$. It follows that
$\mathcal{G}(\mathscr{C}){\text-}\dim A\leq
{^{\bot}\mathscr{C}}{\text-}\dim A$. \hfill$\square$

\vspace{0.5cm}

\centerline{\bf 4. Dual Results}

\vspace{0.2cm}

In this section, we introduce the notion of (pre)coresolvmg
subcategories of $\mathscr{A}$. Then we give some criteria for
computing and comparing homological codimensions relative to
different precoresolving subcategories.  The results and their
proofs in this section are completely dual to that in Section 3, so
we only list the results without proofs.

\vspace{0.2cm}

{\bf Definition 4.1.} Let $\mathscr{E}$ and $\mathscr{T}$ be
subcategories of $\mathscr{A}$. Then $\mathscr{T}$ is called {\it
$\mathscr{E}$-precoresolving} in $\mathscr{A}$ if the following
conditions are satisfied.

(1) $\mathscr{T}$ admits an $\mathscr{E}$-coproper cogenerator.

(2) $\mathscr{T}$ is {\it closed under $\mathscr{E}$-coproper
extensions}, that is, for any
$\Hom_{\mathscr{A}}(-,\mathscr{E})$-exact exact sequence $0\to
A_1\to A_2 \to A_3 \to 0$ in $\mathscr{A}$, if both $A_1$ and $A_3$
are objects in $\mathscr{T}$, then $A_2$ is also an object in
$\mathscr{T}$.

An $\mathscr{E}$-precoresolving subcategory $\mathscr{T}$ of
$\mathscr{A}$ is called {\it $\mathscr{E}$-coresolving} if the
following condition is satisfied.

(3) $\mathscr{T}$ is {\it closed under cokernels of
$\mathscr{E}$-coproper monomorphisms}, that is, for any
$\Hom_{\mathscr{A}}(-,\mathscr{E})$-exact exact sequence $0\to
A_1\to A_2 \to A_3 \to 0$ in $\mathscr{A}$, if both $A_1$ and $A_2$
are objects in $\mathscr{T}$, then $A_3$ is also an object in
$\mathscr{T}$.

\vspace{0.2cm}

The following list shows that the class of
$\mathscr{E}$-(pre)coresolving subcategories is rather large.

\vspace{0.2cm}

{\bf Example 4.2.} (1) Let $\mathscr{A}$ admit enough injective
objects and $\mathscr{E}$ the subcategory of $\mathscr{A}$
consisting of injective objects. Then a subcategory of $\mathscr{A}$
closed under $\mathscr{E}$-coproper extensions is just a subcategory
of $\mathscr{A}$ closed under extensions. Furthermore, if
$\mathscr{C}=\mathscr{E}$ in the above definition, then an
$\mathscr{E}$-precoresolving subcategory is just a subcategory which
contains all injective objects and is closed under extensions, and
an $\mathscr{E}$-coresolving subcategory is just an injectively
coresolving subcategory in the sense of [H2].

(2) Let $\mathscr{C}$ be a subcategory of $\mathscr{A}$ with
$\mathscr{C}\bot\mathscr{C}$. Then by [SSW, Corollary 4.5], the
Gorenstein subcategory $\mathcal{G}(\mathscr{C})$ of $\mathscr{A}$
is a $\mathscr{C}$-precoresolving subcategory of $\mathscr{A}$ with
a $\mathscr{C}$-coproper cogenerator $\mathscr{C}$; furthermore, if
$\mathscr{C}$ is closed under cokernels of monomorphisms, then
$\mathcal{G}(\mathscr{C})$ is a $\mathscr{C}$-coresolving
subcategory of $\mathscr{A}$ by [SSW, Theorem 4.12(b)].

(3) Let $R$ be a ring and $\mathcal{I}(\Mod R)$ the subcategory of
$\Mod R$ consisting of injective modules. If
$(\mathscr{X},\mathscr{Y})$ is a cotorsion pair in $\Mod R$, then
$\mathscr{Y}$ is an $\mathcal{I}(\Mod R)$-precoresolving subcategory
of $\Mod R$ with an $\mathcal{I}(\Mod R)$-coproper cogenerator
$\mathcal{I}(\Mod R)$ ([EJ2]).

(4) Let $R$ be a ring. Recall that a module $E$ in $\Mod R$ is
called {\it FP-injective} if $\Ext_R^1(M, E)=0$ for any finitely
presented left $R$-module $M$. FP-injective modules are also known
as {\it absolutely pure modules}. We use $\mathcal{FI}(\Mod R)$ to
denote the subcategory of $\Mod R$ consisting of FP-injective
modules. Then by [Hu, Lemma 3.1 and Theorem 3.4], it is not
difficult to see that the subcategory of $\Mod R$ consisting of
Gorenstein FP-injective modules (see [MD] or Section 5 below for the
definition) is an $\mathcal{FI}(\Mod R)$-coresolving subcategory of
$\Mod R$ with an $\mathcal{FI}(\Mod R)$-coproper cogenerator
$\mathcal{I}(\Mod R)$.

(5) Let $R$ be a ring. We denote by $\cores
\widetilde{\mathcal{P}(\Mod R)}=\{M \in \Mod R\mid$ there exists a
$\Hom_R(-, \mathcal{P}(\Mod R))$-exact exact sequence $0\to M \to
P^0 \to P^1 \to \cdots \to P^i \to \cdots$ in $\Mod R$ with all
$P^i$ projective$\}$. Then by [Hu, Lemma 3.1 and Theorem 3.8], it is
easy to see that $\cores \widetilde{\mathcal{P}(\Mod R)}$ is a
$\mathcal{P}(\Mod R)$-coresolving subcategory of $\Mod R$ with a
$\mathcal{P}(\Mod R)$-coproper cogenerator $\mathcal{P}(\Mod R)$.
Let $R$ be a left and right noetherian ring. Then by [AB, Theorem
2.17] and [Hu, Lemma 3.1], the subcategory of $\mod R$ consisting of
$\infty$-torsionfree modules (see [HuH] or Section 5 below for the
definition) is a $\mathcal{P}(\mod R)$-coresolving subcategory of
$\mod R$ with a $\mathcal{P}(\mod R)$-coproper cogenerator
$\mathcal{P}(\mod R)$.

\vspace{0.2cm}

Unless stated otherwise, in the rest of this section, we fix a
subcategory $\mathscr{E}$ of $\mathscr{A}$ and an
$\mathscr{E}$-precoresolving subcategory $\mathscr{T}$ of
$\mathscr{A}$ admitting an $\mathscr{E}$-coproper cogenerator
$\mathscr{C}$. We will give some criteria for computing the
$\mathscr{T}$-codimension of a given object $A$ in $\mathscr{A}$,
and then compare it with the $\mathscr{C}$-codimension of $A$.

The following two propositions play a crucial role in this section.

\vspace{0.2cm}

{\bf Proposition 4.3.} {\it Let
$$0\to M \to T^0\to T^1 \to A \to 0\eqno{(4.1)}$$
be an exact sequence in $\mathscr{A}$ with both $T^0$ and $T^1$
objects in $\mathscr{T}$. Then we have

(1) There exists an exact sequence:
$$0\to M\to C \to T\to A\to 0\eqno{(4.2)}$$ in
$\mathscr{A}$ with $T$ an object in $\mathscr{T}$ and $C$ an object
in $\mathscr{C}$.

(2) If (4.1) is $\Hom_{\mathscr{A}}(-,X)$-exact for some object $X$
in $\mathscr{A}$, then so is (4.2).}

\vspace{0.2cm}

As an application of Proposition 4.3, we get the following

\vspace{0.2cm}

{\bf Proposition 4.4.} {\it Let $n\geq 1$ and
$$0\to M\to T^0\to T^1\to \cdots \to T^{n-1}\to A\to 0$$
be an exact sequence in $\mathscr{A}$ with all $T^i$ objects in
$\mathscr{T}$. Then there exist an exact sequence:
$$0\to M\to C^0\to C^1\to \cdots \to C^{n-1}\to B\to 0$$
and a $\Hom_{\mathscr{A}}(-,\mathscr{E})$-exact exact sequence:
$$0\to A\to B\to T\to 0$$ in $\mathscr{A}$ with
all $C^i$ objects in $\mathscr{C}$ and $T$ an object in
$\mathscr{T}$. In particular, an object in $\mathscr{A}$ is an
$n{\text-}\mathscr{C}$-syzygy if and only if it is an
$n{\text-}\mathscr{T}$-syzygy.}

\vspace{0.2cm}

The following corollary is an immediate consequence of Proposition
4.4.

\vspace{0.2cm}

{\bf Corollary 4.5.}  {\it Let $A$ be an object in $\mathscr{A}$ and
$n\geq 0$. If $\mathscr{T}{\text-}\dim A=n$. Then there exists a
$\Hom_{\mathscr{A}}(-,\mathscr{E})$-exact exact sequence $0\to A \to
B \to T \to 0$ in $\mathscr{A}$ with $\mathscr{C}{\text-}\dim B\leq
n$ and $T$ an object in $\mathscr{T}$.}

\vspace{0.2cm}

We give a criterion for computing the $\mathscr{T}$-codimension of
an object in $\mathscr{A}$ as follows.

\vspace{0.2cm}

{\bf Theorem 4.6.} {\it The following statements are equivalent for
any object $M$ in $\mathscr{A}$ and $n \geq 0$.

(1) $\mathscr{T}{\text-}\codim M\leq n$.

(2) There exists an exact sequence: $$0\to M\to C^0\to C^1\to \cdots
\to C^{n-1}\to K^n\to 0$$ in $\mathscr{A}$ with all objects $C^i$ in
$\mathscr{C}$ and $K^n$ an object in $\mathscr{T}$.}

\vspace{0.2cm}

The following result gives a criterion for computing the
$\mathscr{T}$-dimension of an object in $\mathscr{A}$.

\vspace{0.2cm}

{\bf Theorem 4.7.} {\it The following statements are equivalent for
any object $A$ in $\mathscr{A}$ and $n \geq 0$.

(1) $\mathscr{T}{\text-}\dim A\leq n$.

(2) There exists an exact sequence: $$0\to C_n\to C_{n-1}\to
\cdots\to C_1\to K_0 \to A \to 0$$ in $\mathscr{A}$ with $K_0$ an
object in $\mathscr{T}$ and all $C_i$ objects in $\mathscr{C}$, that
is, there exists an exact sequence: $$0\to H\to T\to A \to 0$$ in
$\mathscr{A}$ with $T$ an object in $\mathscr{T}$ and
$\mathscr{C}{\text-}\dim H\leq n-1$.}

\vspace{0.2cm}

The following result gives a sufficient condition such that the
$n{\text-}\mathscr{C}$-cosyzygy of an object in $\mathscr{A}$ with
$\mathscr{T}$-codimension at most $n$ is in $\mathscr{T}$.

\vspace{0.2cm}

{\bf Theorem 4.8.} {\it Let $\mathscr{T}$ be closed under cokernels
of ($\mathscr{E}$-coproper) monomorphisms and
$\mathscr{T}\subseteq{^{\bot}\mathscr{C}}$, and let $n\geq 0$. Then
for any object $M$ in $\mathscr{A}$ with $\mathscr{T}{\text-}\codim
M\leq n$ we have

(1) For any ($\Hom_{\mathscr{A}}(-,\mathscr{E})$-exact) exact
sequence $0\to M \to C^0\to C^1\to \cdots \to C^{n-1}\to K^n \to 0$
in $\mathscr{A}$ with all $C^i$ objects in $\mathscr{C}$, $K^n$ is
an object in $\mathscr{T}$.

(2) If $0\to M\to C \to K \to 0$ is a
($\Hom_{\mathscr{A}}(-,\mathscr{E})$-exact) exact sequence in
$\mathscr{A}$ with $C$ an object in $\mathscr{C}$, then
$\mathscr{T}{\text-}\codim K\leq n-1$.}

\vspace{0.2cm}

%The following corollary is an immediate consequence of Theorem 4.8.

%\vspace{0.2cm}

%{\bf Corollary 4.9.} {\it Let $\mathscr{T}$ admit an injective
%generator and be closed under extensions and cokernels of
%monomorphisms. Then for any object $M$ in $\mathscr{A}$ and $n\geq
%0$ we have

%(1) $\mathscr{T}{\text-}\dim M\leq n$ if and only if for any exact
%sequence $0\to M \to C^0\to C^1\to \cdots \to C^{n-1}\to K^n \to 0$
%in $\mathscr{A}$ with all $C^i$ injective objects in $\mathscr{T}$,
%$K^n$ is an object in $\mathscr{T}$.

%(2) If $0\to M\to I \to K \to 0$ is an exact sequence in
%$\mathscr{A}$ with $\mathscr{T}{\text-}\codim M\leq n$ and $I$ an
%injective object in $\mathscr{T}$, then $\mathscr{T}{\text-}\codim
%K\leq n-1$.}

%\vspace{0.2cm}

We use $\mathscr{T}{\text-}\codim^{\leq n}$ to denote the
subcategory of $\mathscr{A}$ consisting of objects with
$\mathscr{T}$-codimension at most $n$.

\vspace{0.2cm}

{\bf Corollary 4.9.} {\it Let $\mathscr{T}$ be a
$\mathscr{C}$-coresolving subcategory of $\mathscr{A}$ with a
$\mathscr{C}$-coproper cogenerator $\mathscr{C}$ and
$\mathscr{T}\subseteq{^{\bot}\mathscr{C}}$. If $\mathscr{T}$ is
closed under direct summands, then so is
$\mathscr{T}{\text-}\codim^{\leq n}$ for any $n\geq 0$.}

\vspace{0.2cm}

The following result gives some sufficient conditions such that the
$\mathscr{T}$-codimension and the $\mathscr{C}$-codimension of an
object in $\mathscr{A}$ are identical.

\vspace{0.2cm}

{\bf Theorem 4.10.} {\it Let
$\mathscr{T}\subseteq{\mathscr{C}^{\bot}}\bigcap
{^{\bot}\mathscr{C}}$ and $\mathscr{C}$ be closed under direct
summands. Then for an object $M$ in $\mathscr{A}$,
$\mathscr{T}{\text-}\codim M=\mathscr{C}{\text-}\codim M$ if one of
the following conditions is satisfied.

(1) $\mathscr{C}{\text-}\codim M<\infty$, $\mathscr{E}=\mathscr{C}$
and $\mathscr{T}$ is closed under cokernels of
$\mathscr{C}$-coproper monomorphisms.

(2) $\mathscr{C}{\text-}\codim M<\infty$, $\mathscr{E}=\mathscr{C}$
and $\mathscr{C}{\text-}\dim^{<\infty}$ is closed under direct
summands.

(3) $M$ is an object in $^{\bot}\mathscr{T}$ and $\mathscr{C}$ is a
generator for $\mathscr{T}$.}

%\vspace{0.2cm}

%The following result gives another sufficient condition such that
%the $\mathscr{T}$-codimension and the $\mathscr{C}$-codimension of
%an object in $\mathscr{A}$ are identical.

%\vspace{0.2cm}

%{\bf Theorem 4.12.} {\it Let $\mathscr{C}$ be a generator for
%$\mathscr{T}$ such that $\mathscr{T}\subseteq
%{\mathscr{C}^{\bot}}\bigcap {^{\bot}\mathscr{C}}$ and $\mathscr{C}$
%is closed under direct summands. Then for any object $M$ in
%${^{\bot}\mathscr{T}}$, $\mathscr{T}{\text-}\codim
%M=\mathscr{C}{\text-}\codim M$.}

\vspace{0.2cm}

The following result gives a sufficient condition such that the
$\mathscr{T}$-dimension and the $\mathscr{C}$-dimension of an object
in $\mathscr{A}$ are identical.

\vspace{0.2cm}

{\bf Theorem 4.11.} {\it Let $\mathscr{D}$ be a subcategory of
$\mathscr{A}$ such that
$\mathscr{T}\subseteq{\mathscr{C}^{\bot}}\bigcap{^{\bot}\mathscr{D}}$,
and let $\mathscr{C}{\text-}\dim^{\leq n}$ be closed under direct
summands for any $n\geq 0$. If $A$ is an object in $\mathscr{A}$
with $\mathscr{D}{\text-}\codim A<\infty$, then
$\mathscr{T}{\text-}\dim A=\mathscr{C}{\text-}\dim A$.}

\vspace{0.2cm}

In the following, we fix a subcategory $\mathscr{C}$ of
$\mathscr{A}$.

The following two corollaries give some sufficient conditions such
that the $\mathcal{G}(\mathscr{C})$-codimension and the
$\mathscr{C}$-codimension of an object in $\mathscr{A}$ are
identical.

\vspace{0.2cm}

{\bf Corollary 4.12.} {\it Let $\mathscr{C}\bot\mathscr{C}$ and let
$\mathscr{C}$ be closed under direct summands. Then for any object
$M$ in $^{\bot}{\mathcal{G}(\mathscr{C})}$,
$\mathcal{G}(\mathscr{C}){\text-}\codim M=\mathscr{C}{\text-}\codim
M$.}

\vspace{0.2cm}

{\bf Corollary 4.13.} {\it Let $\mathscr{C}\bot\mathscr{C}$ and let
$\mathscr{C}$ be closed under direct summands, and let $\mathscr{D}$
be a subcategory of $^{\bot}{\mathcal{G}(\mathscr{C})}$. Then for
any object $M$ in $^{\bot}\mathscr{C}$ with $\mathscr{D}{\text-}\dim
M<\infty$, $\mathcal{G}(\mathscr{C}){\text-}\codim
M=\mathscr{C}{\text-}\codim M$.}

\vspace{0.2cm}

The following result gives a sufficient condition such that the
$\mathcal{G}(\mathscr{C})$-codimension and the
$\mathscr{C}^{\bot}$-codimension of an object in $\mathscr{A}$ are
identical.

\vspace{0.2cm}

{\bf Theorem 4.14.} {\it Let $\mathscr{C}\bot\mathscr{C}$ and let
$M$ an object in $\mathscr{A}$ with
$\mathcal{G}(\mathscr{C}){\text-}\codim M<\infty$. Then
$\mathcal{G}(\mathscr{C}){\text-}\codim
M={\mathscr{C}^{\bot}}{\text-}\codim M$.}

\vspace{0.5cm}

\centerline{\bf 5. Applications and Questions}

\vspace{0.2cm}

In this section, we will apply the results in Sections 3 and 4 to
special subcategories and in particular to module categories.
Finally we propose some open questions and conjectures concerning
the obtained results.

\vspace{0.2cm}

\noindent{\it 5.1. Special subcategories}

\vspace{0.2cm}

We define $\res\widetilde{\mathscr{C}}=\{A$ is an object in
$\mathscr{A}\mid$ there exists a
$\Hom_{\mathscr{A}}(\mathscr{C},-)$-exact exact sequence $\cdots \to
C_i \to \cdots \to C_1 \to C_0 \to A\to 0$ in $\mathscr{A}$ with all
$C_i$ objects in $\mathscr{C}\}$. Dually, we define
$\cores\widetilde{\mathscr{C}}=\{M$ is an object in
$\mathscr{A}\mid$ there exists a
$\Hom_{\mathscr{A}}(-,\mathscr{C})$-exact exact sequence $0 \to M\to
C^0 \to C^1\to \cdots \to C^i \to \cdots$ in $\mathscr{A}$ with all
$C^i$ objects in $\mathscr{C}\}$ (see [SSW]).

We have the following

\vspace{0.2cm}

{\bf Fact 5.1.} (1) Note that $\mathscr{C}$ is a
$\mathscr{C}$-proper generator for $\res\widetilde{\mathscr{C}}$ and
$\res\widetilde{\mathscr{C}}\bigcap{^{\bot}\mathscr{C}}$. By [Hu,
Lemma 3.1(1)], both $\res\widetilde{\mathscr{C}}$ and
$\res\widetilde{\mathscr{C}}\bigcap{^{\bot}\mathscr{C}}$ are closed
under $\mathscr{C}$-proper extensions. So both
$\res\widetilde{\mathscr{C}}$ and
$\res\widetilde{\mathscr{C}}\bigcap{^{\bot}\mathscr{C}}$ are
$\mathscr{C}$-preresolving. We remark that if $\mathscr{C}$ is a
$\mathscr{C}$-proper generator for $\mathscr{A}$, then
$\res\widetilde{\mathscr{C}}=\mathscr{A}$ and
$\res\widetilde{\mathscr{C}}\bigcap{^{\bot}\mathscr{C}}={^{\bot}\mathscr{C}}$.

(2) If $\mathscr{C}$ is closed under kernels of epimorphisms, then
so are both $\res\widetilde{\mathscr{C}}$ and
$\res\widetilde{\mathscr{C}}\bigcap{^{\bot}\mathscr{C}}$ ([Hu,
Proposition 4.7(1)]).

Dually, we have the following

(3) Note that $\mathscr{C}$ is a $\mathscr{C}$-coproper cogenerator
for $\cores\widetilde{\mathscr{C}}$ and
${\mathscr{C}^{\bot}}\bigcap\cores\widetilde{\mathscr{C}}$. By [Hu,
Lemma 3.1(2)], both $\cores\widetilde{\mathscr{C}}$ and
${\mathscr{C}^{\bot}}\bigcap\cores\widetilde{\mathscr{C}}$ are
closed under $\mathscr{C}$-coproper extensions. So both
$\cores\widetilde{\mathscr{C}}$ and
${\mathscr{C}^{\bot}}\bigcap\cores\widetilde{\mathscr{C}}$ are
$\mathscr{C}$-precoresolving.  We also remark that if $\mathscr{C}$
is a $\mathscr{C}$-coproper cogenerator for $\mathscr{A}$, then
$\cores\widetilde{\mathscr{C}}=\mathscr{A}$ and
${\mathscr{C}^{\bot}}\bigcap\cores\widetilde{\mathscr{C}}={\mathscr{C}^{\bot}}$.

(4) If $\mathscr{C}$ is closed under cokernels of monomorphisms,
then so are both $\cores\widetilde{\mathscr{C}}$ and
${\mathscr{C}^{\bot}}\bigcap\cores\widetilde{\mathscr{C}}$ ([Hu,
Proposition 4.7(2)]).

\vspace{0.2cm}

{\bf Application 5.2.} By Fact 5.1, we can apply the results in
Section 3 in the cases for $\mathscr{T}=\res\widetilde{\mathscr{C}}$
and
$\mathscr{T}=\res\widetilde{\mathscr{C}}\bigcap{^{\bot}\mathscr{C}}$
respectively, and apply the results in Section 4 in the cases for
$\mathscr{T}=\cores\widetilde{\mathscr{C}}$ and
$\mathscr{T}={\mathscr{C}^{\bot}}\bigcap\cores\widetilde{\mathscr{C}}$
respectively. We will not list these consequences in details.

\vspace{0.2cm}

\noindent{\it 5.2. Module categories}

\vspace{0.2cm}

In this subsection, $R$ is a ring and all subcategories of $\Mod R$
are full and additive. For a module $A$ in $\Mod R$, we denote the
projective, injective and flat dimensions of $A$ by $\pd_RA$,
$\id_RA$ and $\fd_RA$ respectively.

We first give the following

\vspace{0.2cm}

{\bf Proposition 5.3.} {\it Let $\mathscr{T}$ be a subcategory of
$\Mod R$.

(1) If $\mathscr{T}$ is closed under extensions and
$\mathcal{P}(\Mod
R)\subseteq\mathscr{T}\subseteq{^{\bot}\mathcal{P}(\Mod R)}$, then
$\pd_RA=\mathscr{T}{\text-}\dim A$ for any $A\in \Mod R$ with
$\pd_RA<\infty$.

(2) If $\mathscr{T}$ is closed under extensions and
$\mathcal{I}(\Mod R)\subseteq\mathscr{T}\subseteq{\mathcal{I}(\Mod
R)^{\bot}}$, then $\id_RA=\mathscr{T}{\text-}\codim A$ for any $A\in
\Mod R$ with $\id_RA<\infty$.}

\vspace{0.2cm}

{\it Proof.} (1) Because
$\mathscr{T}\subseteq{^{\bot}\mathcal{P}(\Mod R)}={\mathcal{P}(\Mod
R)^{\bot}}\bigcap{^{\bot}\mathcal{P}(\Mod R)}$ by assumption, we get
the assertion by Theorem 3.10(2).

(2) It is dual to (1). \hfill$\square$

\vspace{0.2cm}

Let $(\mathscr{X},\mathscr{Y})$ be a cotorsion pair in $\Mod R$.
Then $\mathscr{X}\bigcap\mathscr{Y}$ is called the {\it heart} of
$(\mathscr{X},\mathscr{Y})$. A cotorsion pair
$(\mathscr{X},\mathscr{Y})$ is called {\it hereditary} if
$\mathscr{X}={^{\bot}\mathscr{Y}}$ and
$\mathscr{Y}={\mathscr{X}^{\bot}}$; in this case, $\mathscr{X}$ is
projectively resolving and $\mathscr{Y}$ is injectively coresolving
([GT, Lemma 2.2.10]). By Proposition 5.3 we immediately have the
following

\vspace{0.2cm}

{\bf Corollary 5.4.} {\it Let $(\mathscr{X},\mathscr{Y})$ be a
hereditary cotorsion pair in $\Mod R$ with the heart
$\mathscr{C}(=\mathscr{X}\bigcap\mathscr{Y})$.

(1) If $\mathscr{C}=\mathcal{P}(\Mod R)$, then for any $A\in \Mod R$
with $\pd_RA<\infty$, $\pd_RA=\mathscr{X}{\text -}\dim A$.

(2) If $\mathscr{C}=\mathcal{I}(\Mod R)$, then for any $A\in \Mod R$
with $\id_RA<\infty$, $\id_RA=\mathscr{Y}{\text -}\codim A$.}

\vspace{0.2cm}

Note that a module in $\mathcal{G}(\mathcal{P}(\Mod R))$ (resp.
$\mathcal{G}(\mathcal{I}(\Mod R))$) is just a Gorenstein projective
(resp.  injective) module in $\Mod R$.  So
$\mathcal{G}(\mathcal{P}(\Mod R)){\text-}\dim_RA$ (resp.
$\mathcal{G}(\mathcal{I}(\Mod R)){\text-}\codim_RA$) is just the
Gorenstein projective (resp. injective) dimension of a module $A$ in
$\Mod R$.  We denote the Gorenstein projective (resp. injective)
dimension of a module $A$ in $\Mod R$ by $\Gpd_RA$ (resp.
$\Gid_RA$).

\vspace{0.2cm}

{\bf Definition 5.5.} Let $A$ be a module in $\Mod R$.

(1) ([DLM]) $A$ is called {\it strongly Gorenstein flat} if there
exists a $\Hom_R(-,\mathcal{F}(\Mod R))$-exact exact sequence
$\cdots \to P_1 \to P_0 \to P^0 \to P^1 \to \cdots$ in $\Mod R$ with
all terms projective, such that $A\cong\Im(P_0 \to P^0)$. We use
$\mathcal{SGF}(\Mod R)$ to denote the subcategory of $\Mod R$
consisting of strongly Gorenstein flat modules.

The {\it strongly Gorenstein flat dimension} $\SGfd_RA$ of $A$ is
defined to be $\inf\{n\mid$ there exists an exact sequence $0\to G_n
\to \cdots \to G_1 \to G_0 \to A \to 0$ in $\Mod R$ with all $G_i$
in $\mathcal{SGF}(\Mod R)\}$. Set $\SGfd_RA=\infty$ if no such $n$
exists.

(2) ([MD]) $A$ is called {\it Gorenstein FP-injective} if there
exists a $\Hom_R(\mathcal{FI}(\Mod R),-)$-exact exact sequence
$\cdots \to I_1 \to I_0 \to I^0 \to I^1 \to \cdots$ in $\Mod R$ with
all terms injective, such that $A\cong\Im(I_0 \to I^0)$. We use
$\mathcal{GFI}(\Mod R)$ to denote the subcategory of $\Mod R$
consisting of Gorenstein FP-injective modules.

The {\it Gorenstein FP-injective dimension} $\GFid_RA$ of $A$ is
defined to be $\inf\{n\mid$ there exists an exact sequence $0\to A
\to H^0 \to H^1\to \cdots \to H^n \to 0$ in $\Mod R$ with all $H^i$
in $\mathcal{GFI}(\Mod R)\}$. Set $\GFid_A=\infty$ if no such $n$
exists.

\vspace{0.2cm}

It is trivial that there exist the following inclusions:
{\footnotesize$$\mathcal{P}(\Mod R)\subseteq\mathcal{SGF}(\Mod
R)\subseteq\mathcal{G}(\mathcal{P}(\Mod R))\left\{
\begin{array}{ll}\subseteq \cores\widetilde{\mathcal{P}(\Mod
R)},\\
\subseteq{^{\bot}(\mathcal{P}(\Mod
R))}\supseteq{^{\bot}(\mathcal{F}(\Mod
R))}\supseteq\mathcal{SGF}(\Mod R),\end{array} \right.$$} and
{\footnotesize$$\mathcal{I}(\Mod R)\subseteq\mathcal{GFI}(\Mod
R)\subseteq\mathcal{G}(\mathcal{I}(\Mod R))\left\{
\begin{array}{ll}\subseteq \res\widetilde{\mathcal{I}(\Mod
R)},\\
\subseteq{(\mathcal{I}(\Mod R))^{\bot}}\supseteq{(\mathcal{FI}(\Mod
R))^{\bot}}\supseteq\mathcal{GFI}(\Mod R).\end{array} \right.$$} So
for any module $A$ in $\Mod R$, we have
{\footnotesize$$\pd_RA\geq\SGfd_RA\geq\Gpd_RA\left\{
\begin{array}{ll}\geq \cores\widetilde{\mathcal{P}(\Mod
R)}{\text -}\dim A,\\
\geq{^{\bot}(\mathcal{P}(\Mod R))}{\text -}\dim
A\leq{^{\bot}(\mathcal{F}(\Mod R))}{\text -}\dim
A\leq\SGfd_RA,\end{array} \right.$$} and
{\footnotesize$$\id_RA\geq\GFid_RA\geq\Gid_RA\left\{
\begin{array}{ll}\geq \res\widetilde{\mathcal{I}(\Mod
R)}{\text -}\codim A,\\
\geq{(\mathcal{I}(\Mod R))^{\bot}}{\text -}\codim
A\leq{(\mathcal{FI}(\Mod R))^{\bot}}{\text -}\codim
A\leq\GFid_RA.\end{array} \right.$$}

\vspace{0.2cm}

{\bf Theorem 5.6.} {\it Let $A$ be a module in $\Mod R$.

(1) If $A\in ({\mathcal{SGF}(\Mod R)})^{\bot}$, then
$\pd_RA=\SGfd_RA$.

(2) If $A\in (\mathcal{G}({\mathcal{P}(\Mod R))})^{\bot}$, then
$\pd_RA=\SGfd_RA=\Gpd_RA$.

(3) If $\id_RA<\infty$, then $\pd_RA=\SGfd_RA=\Gpd_RA=
\cores\widetilde{\mathcal{P}(\Mod R)}{\text -}\dim A$.

(4) If $\pd_RA<\infty$, then $\pd_RA=\SGfd_RA=\Gpd_RA=
{^{\bot}({\mathcal{P}(\Mod R)})}{\text -}\dim
A={^{\bot}({\mathcal{F}(\Mod R)})}{\text -}\dim A$.

(5) If $\SGfd_RA<\infty$, then $\SGfd_RA=\Gpd_RA=
{^{\bot}({\mathcal{P}(\Mod R)})}{\text -}\dim
A={^{\bot}({\mathcal{F}(\Mod R)})}$-$\dim A$.

(6) If $\Gpd_RA<\infty$, then $\Gpd_RA={^{\bot}({\mathcal{P}(\Mod
R)})}{\text -}\dim A$.

(7) If $\fd_RA<\infty$, then $\pd_RA=\SGfd_RA$.}

\vspace{0.2cm}

{\it Proof.}  (1) (resp. (2)) It is clear that $\mathcal{P}(\Mod R)$
is both a $\mathcal{P}(\Mod R)$-proper generator and a
$\mathcal{P}(\Mod R)$-coproper cogenerator for $\mathcal{SGF}(\Mod
R)$ (resp. $\mathcal{G}(\mathcal{P}(\Mod R))$). Note that
$\mathcal{SGF}(\Mod R)$ (resp. $\mathcal{G}(\mathcal{P}(\Mod R))$)
is closed under extensions by [Hu, Lemma 3.1] (resp. [H2, Theorem
2.5]). Then by putting $\mathscr{T}=\mathcal{SGF}(\Mod R)$ (resp.
$\mathscr{T}=\mathcal{G}(\mathcal{P}(\Mod R))$) and
$\mathscr{C}=\mathcal{P}(\Mod R)$ in Theorem 3.10(3), we have
$\pd_RA=\SGfd_RA$ (resp. $\pd_RA=\Gpd_RA$) if $A\in
({\mathcal{SGF}(\Mod R)})^{\bot}$ (resp.
$A\in(\mathcal{G}(\mathcal{P}(\Mod R)))^{\bot}$).

(3) It is clear that $\mathcal{P}(\Mod R)$ is a $\mathcal{P}(\Mod
R)$-coproper cogenerator for $\cores\widetilde{\mathcal{P}(\Mod
R)}$. Note that $\cores\widetilde{\mathcal{P}(\Mod R)}$ is closed
under $\mathcal{P}(\Mod R)$-coproper extensions by [Hu, Lemma 3.1].
Then by putting $\mathscr{T}=\cores\widetilde{\mathcal{P}(\Mod R)}$,
$\mathscr{C}=\mathcal{P}(\Mod R)$ and $\mathscr{D}=\mathcal{I}(\Mod
R)$ in Theorem 4.11, we have
$\pd_RA=\cores\widetilde{\mathcal{P}(\Mod R)}{\text -}\dim A$ if
$\id_RA<\infty$.

(4) By Proposition 5.3(1), we have
$\pd_RA={^{\bot}({\mathcal{P}(\Mod R)})}{\text -}\dim A$ if
$\pd_RA<\infty$.

(5) Note that $\mathcal{SGF}(\Mod R)$ is closed under extensions
(see the proof of (1)). Let $\SGfd_RA=n<\infty$. Then by Theorem
3.6, there exists an exact sequence:
$$0\to G_n\to P_{n-1}\to P_{n-2}\to \cdots \to P_0\to A\to 0$$
in $\Mod R$ with all $P_i$ projective and $G_n$ strongly Gorenstein
flat. Put $K_i=\Im(P_i\to P_{i-1})$ for any $1\leq i\leq n-1$.
Suppose ${^{\bot}({\mathcal{P}(\Mod R)})}{\text -}\dim A=m<\infty$.
It suffices to show $m\geq n$. If $m<n$, then $K_m
\in{^{\bot}({\mathcal{P}(\Mod R)})}$ and $K_{n-1}
\in{^{\bot}({\mathcal{P}(\Mod R)})}$ by Theorem 3.8(1). Because
there exists a $\Hom_R(-,\mathcal{F}(\Mod R))$-exact exact sequence
$0 \to G_n\to P \to G\to 0$ in $\Mod R$ with $P$ projective and $G$
strongly Gorenstein flat, we have the following push-out diagram:
$$\xymatrix{
& 0 \ar[d] &0 \ar[d] & & \\
0 \ar[r] & G_n \ar[d] \ar[r] & P_{n-1}
\ar[d] \ar[r] & K_{n-1} \ar@{=}[d] \ar[r]  &0  \\
0 \ar[r] & P \ar[d] \ar[r] & G^{'}
\ar[d] \ar[r] & K_{n-1} \ar[r]  &0  \\
& G\ar[d]\ar@{=}[r] & G \ar[d] &   &    \\
&  0 & 0  &  & }$$ By using an argument similar to that in the proof
of [H2, Theorem 2.5], we get that $\mathcal{SGF}(\Mod R)$ is closed
under direct summands. Because both the middle column and the middle
row in the above diagram split, $G^{'}\cong P_{n-1}\oplus G$ is
strongly Gorenstein flat and $K_{n-1}$ is isomorphic to a direct
summand of $G^{'}$, which implies that $K_{n-1}$ is strongly
Gorenstein flat and $\SGfd_RA\leq n-1$. It is a contradiction.

(6) We get the assertion by putting $\mathscr{C}=\mathcal{P}(\Mod
R)$ in Theorem 3.14 or $\mathscr{C}=\mathcal{P}(\Mod R)$ and
$\mathscr{T}=\mathcal{G}(\mathcal{P}(\Mod R))$ in Theorem 3.10.

(7) By the definition of strongly Gorenstein flat modules, it is
easy to see that $A\in ({\mathcal{SGF}(\Mod R)})^{\bot}$ if
$\fd_RA<\infty$. Then the assertion follows from (1).
\hfill$\square$

\vspace{0.2cm}

{\it Remark 5.7.} Theorem 5.6(2) is [Z, Theorem 2.3]. Theorem 5.6(3)
generalizes [H1, Theorem 2.2] which states that for a module $A$ in
$\Mod R$, if $\id_RA<\infty$, then $\pd_RA=\Gpd_RA$. Notice that a
module in $\Mod R$ with finite injective dimension is in
$(\mathcal{G}({\mathcal{P}(\Mod R))})^{\bot}$, so we may also get
[H1, Theorem 2.2] by Theorem 5.6(2) ([Z, Corollary 2.5]). Theorem
5.6(6) is well known ([H2, Theorem 2.20]).

\vspace{0.2cm}

Let $A$ be a module in $\Mod R$. Recall that $A$ is called {\it
Gorenstein flat} if there exists an exact sequence:
$$\cdots \to F_1 \to F_0 \to F^0 \to F^1 \to \cdots$$
in $\Mod R$ with all terms flat, such that $A\cong\Im(F_0 \to F^0)$
and the sequence remains still exact after applying the functor
$I\bigotimes _R-$ for any injective right $R$-module $I$. The {\it
Gorenstein flat dimension} of $A$, denoted by $\Gfd_RA$, is defined
as $\inf\{ n\mid$ there exists an exact sequence $0 \to H_n \to
\cdots \to H_1 \to H_0 \to A \to 0$ with all $H_i$ Gorenstein
flat$\}$. Set $\Gfd_RA=\infty$ if no such $n$ exists. ([EJT, H2]).

The following is an open question: whether is every Gorenstein
projective module over any ring Gorenstein flat? Holm proved in [H2,
Proposition 2.4] that if $R$ is a right coherent ring with finite
left finitistic projective dimension, then every Gorenstein
projective module in $\Mod R$ is Gorenstein flat. As an immediate
consequence of Theorem 5.6, we have the following

\vspace{0.2cm}

{\bf Corollary 5.8.} {\it Let $R$ be a right coherent ring and $A$ a
module in $\Mod R$. Then $\Gpd_RA\geq\Gfd_RA$ if either of the
following conditions is satisfied:

(1) $A\in (\mathcal{G}({\mathcal{P}(\Mod R))})^{\bot}$ (in
particular, if $\pd_RA<\infty$ or $\id_RA<\infty$),

(2) $\SGfd_RA<\infty$.}

\vspace{0.2cm}

{\it Proof.} Let $R$ be a right coherent ring and $A$ a module in
$\Mod R$. Then $\SGfd_RA\geq\Gfd_RA$ by [DLM, Proposition 2.3]. So
the assertions follow from Theorem 5.6(2)(5) respectively.
\hfill$\square$

\vspace{0.2cm} Recall from [B1] that $R$ is called {\it left
GF-closed} if the subcategory of $\Mod R$ consisting of Gorenstein
flat modules is closed under extensions.

\vspace{0.2cm}

{\bf Corollary 5.9.} {\it Let $A$ be a module in $\Mod R$ and $n$ a
non-negative integer.

(1) ([CFH, Lemma 2.17]) If $\Gpd_RA=n$, then there exists an exact
sequence $0\to A\to B \to T \to 0$ in $\Mod R$ with $\pd_BA=n$ and
$T$ Gorenstein projective.

(2) If $\SGfd_RA=n$, then there exists an exact sequence $0\to A\to
B \to T \to 0$ in $\Mod R$ with $\pd_BA=n$ and $T$ strongly
Gorenstein flat.

(3) If $R$ is left GF-closed and $\Gfd_RA=n$, then there exists an
exact sequence $0\to A\to B \to T \to 0$ in $\Mod R$ with $\fd_RB=n$
and $T$ Gorenstein flat.}

\vspace{0.2cm}

{\it Proof.} (1) (resp. (2)) Let $\Gpd_RA=n$ (resp. $\SGfd_RA=n$).
By Corollary 4.5, there exists an exact sequence $0\to A\to B \to T
\to 0$ in $\Mod R$ with $\pd_RB\leq n$ and $T$ Gorenstein projective
(resp. strongly Gorenstein flat). Then by [H2, Theorem 2.5] (resp.
Example 3.2(4)) and Proposition 2.3, we have $\Gpd_RB\geq \Gpd_RA=n$
(resp. $\SGfd_RB\geq \SGfd_RA=n$). So $\pd_RB=n$ by Theorem 5.6(4).

(3) Let $R$ be left GF-closed and $\Gfd_RA=n$. By Corollary 4.5,
there exists an exact sequence $0\to A\to B \to T \to 0$ in $\Mod R$
with $\fd_RB\leq n$ and $T$ Gorenstein flat. Then by [B1, Theorem
2.3] and Proposition 2.3, we have $\Gfd_RB\geq \Gfd_RA=n$. So
$\fd_RB=n$ by [B2, Theorem 2.2]. \hfill$\square$

\vspace{0.2cm}

Recall that the {\it FP-injective dimension} $\FP{\text -}\id_RA$ of
$A$ in $\Mod R$ is defined to be $\inf\{n\mid$ there exists an exact
sequence $0\to A \to E^0 \to E^1\to \cdots \to E^n \to 0$ in $\Mod
R$ with all $E^i$ in $\mathcal{FI}(\Mod R)\}$. Set $\FP{\text
-}\id_RA=\infty$ if no such $n$ exists.

The following result is the dual of Theorem 5.6.

\vspace{0.2cm}

{\bf Theorem 5.10.} {\it Let $A$ be a module in $\Mod R$.

(1) If $A\in {^{\bot}({\mathcal{GFI}(\Mod R)})}$, then
$\id_RA=\GFid_RA$.

(2) If $A\in {^{\bot}(\mathcal{G}({\mathcal{I}(\Mod R))})}$, then
$\id_RA=\GFid_RA=\Gid_RA$.

(3) If $\pd_RA<\infty$, then $\id_RA=\GFid_RA=\Gid_RA=
\res\widetilde{\mathcal{I}(\Mod R)}{\text -}\codim A$.

(4) If $\id_RA<\infty$, then $\id_RA=\GFid_RA=\Gid_RA=
{({\mathcal{I}(\Mod R)})^{\bot}}{\text -}\codim
A={({\mathcal{FI}(\Mod R)})^{\bot}}{\text -}\codim A$.

(5) If $\GFid_RA<\infty$, then $\GFid_RA=\Gid_RA=
{({\mathcal{I}(\Mod R)})^{\bot}}{\text -}\codim
A={({\mathcal{FI}(\Mod R)})^{\bot}}$-$\codim A$.

(6) If $\Gid_RA<\infty$, then $\Gid_RA={({\mathcal{I}(\Mod
R)}){^{\bot}}}{\text -}\codim A$.

(7) If $\FP{\text -}\id_RA<\infty$, then $\id_RA=\GFid_RA$.}

\vspace{0.2cm}

{\it Remark 5.11.} Theorem 5.10(3) generalizes [H1, Theorem 2.1]
which states that for a module $A$ in $\Mod R$, if $\pd_RA<\infty$,
then $\id_RA=\Gid_RA$. Notice that a module in $\Mod R$ with finite
projective dimension is in $^{\bot}(\mathcal{G}({\mathcal{I}(\Mod
R))})$, so we may also get [H1, Theorem 2.1] by Theorem 5.10(2).
Theorem 5.10(6) is well known ([H2, Theorem 2.22]).

\vspace{0.2cm}

\noindent{\it 5.3. Questions}

\vspace{0.2cm}

In view of the assertions (3) and (4) in Theorem 5.6, it is natural
to ask the following

\vspace{0.2cm}

{\bf Question 5.12.} {\it If $A$ is a module in $\Mod R$ with
$\id_RA<\infty$, does then $\pd_RA={^{\bot}(\mathcal{P}(\Mod
R))}{\text -}\dim A$ hold?}

\vspace{0.2cm}

{\bf Question 5.13.} {\it If $A$ is a module in $\Mod R$ with
$\pd_RA<\infty$, does then $\pd_RA=\cores
\widetilde{\mathcal{P}(\Mod R)}{\text -}\dim A$ hold?}

\vspace{0.2cm}

From now on, $R$ is a left and right Noetherian ring (unless stated
otherwise). We denote by ${^{\bot}{_RR}}=\{M\in \mod R \mid \Ext
_R^i(_RM, {_RR})=0$ for any $i\geq 1\}$ (resp. ${^{\bot}R_R}=\{N\in
\mod R^{op}\mid \Ext^i_{R^{op}}(N_R, {R_R})=0$ for any $i\geq 1\}$).

For any module $A$ in $\mod R$, there exists a projective
presentation:
$$P_1\buildrel {f} \over \longrightarrow P_0 \to A \to 0$$
of $A$ in $\mod R$ (note: if $R$ is an artinian algebra, then this
projective presentation of $A$ is chosen to be the minimal one).
Then we get an exact sequence:
$$0\to A^* \to P_0^*\buildrel {f^*} \over \longrightarrow P_1^* \to \Tr A \to 0$$
in $\mod R^{op}$, where $(-)^*=\Hom(-,R)$ and $\Tr A=\Coker f^*$ is
the {\it transpose} of $A$. Auslander and Bridger generalized the
notions of finitely generated projective modules and the projective
dimension of finitely generated modules as follows. A module $A$ in
$\mod R$ is said to {\it have Gorenstein dimension zero} if
$A\in{^{\bot}{_RR}}$ and $\Tr A \in{^{\bot}R_R}$ ([AB]). It is well
known that over a left and right Noetherian ring, a finitely
generated module is Gorenstein projective if and only if it has
Gorenstein dimension zero ([EJ2, Proposition 10.2.6]).

Let $A$ be a module in $\mod R$. Recall from [HuH] that $A$ is
called {\it $\infty$-torsionfree} if $\Tr A \in{^{\bot}R_R}$. We use
$\mathcal{T}(\mod R)$ to denote the subcategory of $\mod R$
consisting of $\infty$-torsionfree modules. The {\it torsionfree
dimension} of $A$, denoted by $\mathcal{T}{\text -}\dim_RA$, is
defined as $\inf\{n\mid$ there exists an exact sequence $0\to X_n
\to \cdots \to X_1 \to X_0 \to A \to 0$ in $\mod R$ with all $X_i$
in $\mathcal{T}(\mod R)\}$. Set $\mathcal{T}{\text -}\dim_RA=\infty$
if no such $n$ exists. By [AB, Theorem 2.17], a module is in $\cores
\widetilde{\mathcal{P}(\mod R)}$ if and only if it is in
$\mathcal{T}(\mod R)$. So $\cores \widetilde{\mathcal{P}(\mod
R)}{\text -}\dim A=\mathcal{T}{\text -}\dim_RA$ for any module $A$
in $\mod R$.

By Example 4.2(5) and Corollary 4.5, we immediately have the
following

\vspace{0.2cm}

{\bf Corollary 5.14.} ([HuH, Corollary 3.5]) {\it Let $A$ be a
module in $\mod R$ with $\mathcal{T}{\text -}\dim_RA=n$. Then there
exists an exact sequence $0\to A\to B \to T \to 0$ in $\Mod R$ with
$\pd_RB\leq n$ and $T$ $\infty$-torsionfree.}

\vspace{0.2cm}

The following result is analogous to Theorem 5.6(3)(4).

\vspace{0.2cm}

{\bf Theorem 5.15.} {\it Let $R$ be a left and right Noetherian ring
and $A$ a module in $\mod R$.

(1) If $\id_RA<\infty$, then $\pd_RA=\Gpd_RA=\mathcal{T}{\text
-}\dim_RA$.

(2) If $\pd_RA<\infty$, then $\pd_RA=\Gpd_RA={^{\bot}{_RR}}{\text
-}\dim A$.}

\vspace{0.2cm}

In view of the assertions in Theorem 5.15, it is natural to ask the
following questions, which are finitely generated versions of
Questions 5.12 and 5.13 respectively.

\vspace{0.2cm}

{\bf Question 5.16.} {\it If $A$ is a module in $\mod R$ with
$\id_RA<\infty$, does then $\pd_RA={^{\bot}{_RR}}{\text -}\dim A$
hold? (equivalently, does then $\Gpd_RA={^{\bot}{_RR}}{\text -}\dim
A$ hold?)}

\vspace{0.2cm}

{\bf Question 5.17.} {\it If $A$ is a module in $\mod R$ with
$\pd_RA<\infty$, does then $\pd_RA=\mathcal{T}{\text -}\dim_RA$
hold? (equivalently, does then $\Gpd_RA=\mathcal{T}{\text -}\dim_RA$
hold?)}

\vspace{0.2cm}

Let $R$ be an artinian algebra and $C(R)$ the center of $R$, and let
$J$ be the Jacobson radical of $C(R)$ and $I(C(R)/J)$ the injective
envelope of $C(R)/J$. Then the Matlis duality
$\mathbb{D}(-)=\Hom_{C(R)}(-,I(C(R)/J))$ between $\mod R$ and $\mod
R^{op}$ induces a duality between projective (resp. injective)
modules in $\mod R$ and injective (resp. projective) modules in
$\mod R^{op}$. As a special case of Question 5.16, we propose the
following

\vspace{0.2cm}

{\bf Conjecture 5.18.} {\it Let $R$ be an artinian algebra.

(1) A module $A$ in $\mod R$ is projective if $A$ is injective and
$A\in {^{\bot}{_RR}}$.

(2) $R$ is selfinjective if $\mathbb{D}(R_R)\in {^{\bot}{_RR}}$.}

\vspace{0.2cm}

The generalized Nakayama conjecture ({\bf GNC} for short) states
that over any artinian algebra $R$, a module $A$ in $\mod R$ is
projective if $\Ext_R^i(A\oplus R, A\oplus R)=0$ for any $i \geq 1$
([AR]). The strong Nakayama conjecture ({\bf SNC} for short) states
that over any artinian algebra $R$, for any $0\neq A$ in $\mod R$
there exists an $i\geq 0$ such that $\Ext_R^i(A,R)\neq 0$ ([CoF]).
These two conjectures remain still open. Observe that an equivalent
version of {\bf GNC} states that over any artinian algebra $R$, for
any simple module $S$ in $\mod R$ there exists $i\geq 0$ such that
$\Ext_R^i(S,R)\neq 0$ ([AR]). So {\bf SNC}$\Rightarrow${\bf GNC}. It
is easy to see that {\bf GNC}$\Rightarrow$Conjecture
5.18(1)$\Rightarrow$Conjecture 5.18(2).

The following result shows that Question 5.17 is closely related to
{\bf SNC}.

\vspace{0.2cm}

{\bf Proposition 5.19.} {\it Let $R$ be an artinian algebra. Then
the following statements are equivalent.

(1) {\bf SNC} holds for $R^{op}$.

(2) A module in $\mod R$ is projective if $A\in \mathcal{T}(\mod R)$
and $\pd_RA\leq 1$.}

\vspace{0.2cm}

{\it Proof.} $(1)\Rightarrow (2)$ Let $A\in\mathcal{T}(\mod R)$ and
$\pd_RA\leq 1$. Then $\Tr A \in{^{\bot}R_R}$ and there exists a
minimal projective presentation:
$$0\to P_1 \to P_0 \to A \to 0$$
in $\mod R$, which induces an exact sequence:
$$0\to A^* \to P_0^*\to P_1^* \to \Tr A \to 0$$
in $\mod R^{op}$. So we get the following commutative diagram with
exact rows:
$$\xymatrix{& 0 \ar[r]\ar[d] & P_1 \ar[r]
\ar[d]^{\cong} & P_0\ar[r] \ar[d]^{\cong} &
A \ar[r] & 0 \\
0\ar[r] & (\Tr A)^* \ar[r] & P_1^{**} \ar[r] & P_0^{**} & & }$$ Thus
$(\Tr A)^*=0$ and $\Ext_R^i(\Tr A, R)=0$ for any $i\geq 0$. Then
$\Tr A=0$ by (1), which implies that $A$ is projective by [ARS,
Chapter IV, Proposition 1.7(b)].

$(2)\Rightarrow (1)$ Let $B$ be a module in $\mod R^{op}$ such that
$\Ext_{R^{op}}^i(B,R)=0$ for any $i\geq 0$. Then $B$ has no non-zero
projective summands. So $B\cong \Tr\Tr B$ by [ARS, Chapter IV,
Proposition 1.7(c)] and hence $\Ext_{R^{op}}^i(\Tr\Tr B,R)=0$ for
any $i\geq 0$. So $\Tr B \in \mathcal{T}(\mod R)$. From a minimal
projective presentation $Q_1\to Q_0 \to B \to 0$ of $B$ in $\mod
R^{op}$, we get an exact sequence:
$$0\to B^* \to Q_0^*\to Q_1^* \to \Tr B \to 0$$
in $\mod R$ with $Q_0^*, Q_1^*$ projective. Because $B^*=0$,
$\pd_R\Tr B\leq 1$. Then $\Tr B$ is projective by (2), which implies
that $B$ is projective. Again because $B^*=0$, $B=0$. Therefore {\bf
SNC} holds for $R^{op}$. \hfill$\square$

\vspace{0.5cm}

{\bf Acknowledgements.} This research was partially supported by the
Specialized Research Fund for the Doctoral Program of Higher
Education (Grant No. 20100091110034), NSFC (Grant No. 11171142), NSF
of Jiangsu Province of China (Grant No. BK2010007) and a Project
Funded by the Priority Academic Program Development of Jiangsu
Higher Education Institutions. The author thanks the referee for the
useful suggestions.


\vspace{0.5cm}

\begin{thebibliography}{101}

\bibitem[A]{A1} M. Auslander, {\it Coherent functors}, in: Proc. of the Conf. on Categorial Algebra,
La Jolla, 1965, Springer-Verlag, Berlin, 1966, pp.189--231.

\bibitem[AB]{A2} M. Auslander and M. Bridger, {\it Stabe Module Theory},
Memoirs Amer. Math. Soc. {\bf 94}, Amer. Math. Soc., Providence, RI,
1969.

\bibitem[AR]{A3} M. Auslander and I. Reiten, {\it On a generalized version of the
Nakayama conjecture}, Proc. Amer. Math. Soc. {\bf 52} (1975),
69--74.

\bibitem[ARS]{A4} M. Auslander, I. Reiten and S.O. Smol$\phi$, Representation Theory
of Artin Algebras, Corrected reprint of the 1995 original, Cambridge
Studies in Adv. Math. {\bf 36}, Cambridge Univ. Press, Cambridge,
1997.

\bibitem[B1]{A5} D. Bennis, {\it Rings over which the class of Gorenstein
flat modules is closed under extensions}, Comm. Algebra {\bf 37}
(2009), 855--868.

\bibitem[B2]{A6} D. Bennis, {\it A note on Gorenstein flat dimension}, Algebra Colloq.
{\bf 18} (2011), 155--162.

\bibitem[C]{A7} L.W. Christensen, Gorenstein Dimension, Lect. Notes in Math. {\bf 1747},
Springer-Verlag, Berlin, 2000.

\bibitem[CFH]{A8} L.W. Christensen, A. Frankild and H. Holm, {\it On Gorenstein projective,
injective and flat dimensions--A functorial description with
applications}, J. Algebra {\bf 302} (2006), 231--279.

\bibitem[CI]{A9} L.W. Christensen and S. Iyengar, {\it Gorenstein
dimension of modules over homomorphisms}, J. Pure Appl. Algebra {\bf
208} (2007), 177--188.

\bibitem[CoF]{A10} R.R. Colby and K.R. Fuller, {\it A note on the Nakayama conjectures},
Tsukaba J. Math. {\bf 14} (1990), 343--352.

\bibitem[DLM]{A11} N.Q. Ding, Y.L. Li and L.X. Mao, {\it Strongly Gorenstein flat
modules}, J. Aust. Math Soc. {\bf 86} (2009), 323--338.

\bibitem[EJ1]{A12} E.E. Enochs and O.M.G. Jenda, {\it Gorenstein injective and
projective modules}, Math. Z. {\bf 220} (1995), 611--633.

\bibitem[EJ2]{A13} E.E. Enochs and O.M.G. Jenda, Relative
Homological Algebra, De Gruyter Exp. in Math. {\bf 30}, Walter de
Gruyter, Berlin, New York, 2000.

\bibitem[EJL]{A14} E.E. Enochs, O.M.G. Jenda and J. A. L\'{o}pez-ramos,
{\it Covers and envelopes by $V$-Gorenstein modules}, Comm. Algebra
{\bf 33} (2005), 4705--4717.

\bibitem[EJT]{A15} E.E. Enochs, O.M.G. Jenda and B. Torrecillas,
{\it Gorenstein flat modules}, J. Nanjing Univ. Math. Biquarterly
{\bf 10} (1993), 1--9.

\bibitem[GD]{A16} Y.X. Geng and N.Q. Ding, {\it $\mathcal{W}$-Gorenstein
moduless}, J. Algebra {\bf 325} (2011), 132--146.

\bibitem[GT]{A17} R. G\"obel and J. Triifaj, Approximations and Endomorphism
Algebras of Modules, De Gruyter Exp. in Math. {\bf 41}, Walter de
Gruyter, Berlin, New York, 2006.

\bibitem[H1]{A18} H. Holm, {\it Rings with finite Gorenstein injective dimension},
Proc. Amer. Math. Soc. {\bf 132} (2003), 1279--1283.

\bibitem[H2]{A19} H. Holm, {\it Gorenstein homological dimensions}, J.
Pure Appl. Algebra {\bf 189} (2004), 167--193.

\bibitem[Hu]{A20} Z.Y. Huang, {\it Proper resolutions and Gorenstein categories},
J. Algebra {\bf 393} (2013), 142--169.

\bibitem[HuH]{A21} C.H. Huang and Z.Y. Huang, {\it Torsionfree dimension of
modules and self-injective dimension of rings}, Osaka J. Math. {\bf
49} (2012), 21--35.

\bibitem[LHX]{A22} Z.F. Liu, Z.Y. Huang and A.M. Xu, {\it Gorenstein projective dimension
relative to a semidualizing bimodule}, Comm. Algebra {\bf 41}
(2013), 1--18.

\bibitem[MD]{A23} L.X. Mao and N.Q. Ding, {\it Gorenstein FP-injective and
Gorenstein flat modules}, J. Algebra Appl. {\bf 7} (2008), 491--606.

\bibitem[SSW]{A24} S. Sather-Wagstaff, T. Sharif and D. White, {\it Stability
of Gorenstein categories}, J. London Math. Soc. {\bf 77} (2008),
481--502.

\bibitem[Z]{A25} P. Zhang, {\it Gorensteinness and Buchweitz theorem}, Preprint 2012.

\end{thebibliography}
\end{document}